\documentclass{vgtc}                          

\graphicspath{{figures/}{pictures/}{images/}{./}} 

\usepackage{times}                     

\usepackage{tabu}                      
\usepackage{booktabs}                  
\usepackage{lipsum}                    
\usepackage{mwe}                       
\usepackage{subfigure}
\usepackage{mathptmx}                  
\usepackage{amsmath}
\usepackage{amsfonts}

\onlineid{0}

\usepackage{enumitem}

\vgtccategory{Research}

\vgtcinsertpkg

\title{Visualizing Local Maxima of the Ohio overdose epidemic with Vineyards}

\author{Nicholas Bermingham \thanks{e-mail: bermingham.11@osu.edu}\\ %
        \scriptsize The Ohio State University %
\and David White\thanks{e-mail: davw710@gmail.com}\\ %
      \scriptsize Queen's University \\ \scriptsize Denison University %
\and Nathan Willey \thanks{e-mail: willey.106@osu.edu}\\ %
     \parbox{1.4in}{\scriptsize \centering The Ohio State University}}

\teaser{
  \centering
  \begin{minipage}{.4\linewidth}
    \centering
    \includegraphics[
        width=\linewidth,
        trim={3cm 1.5cm 1cm 3cm},
        clip
    ]{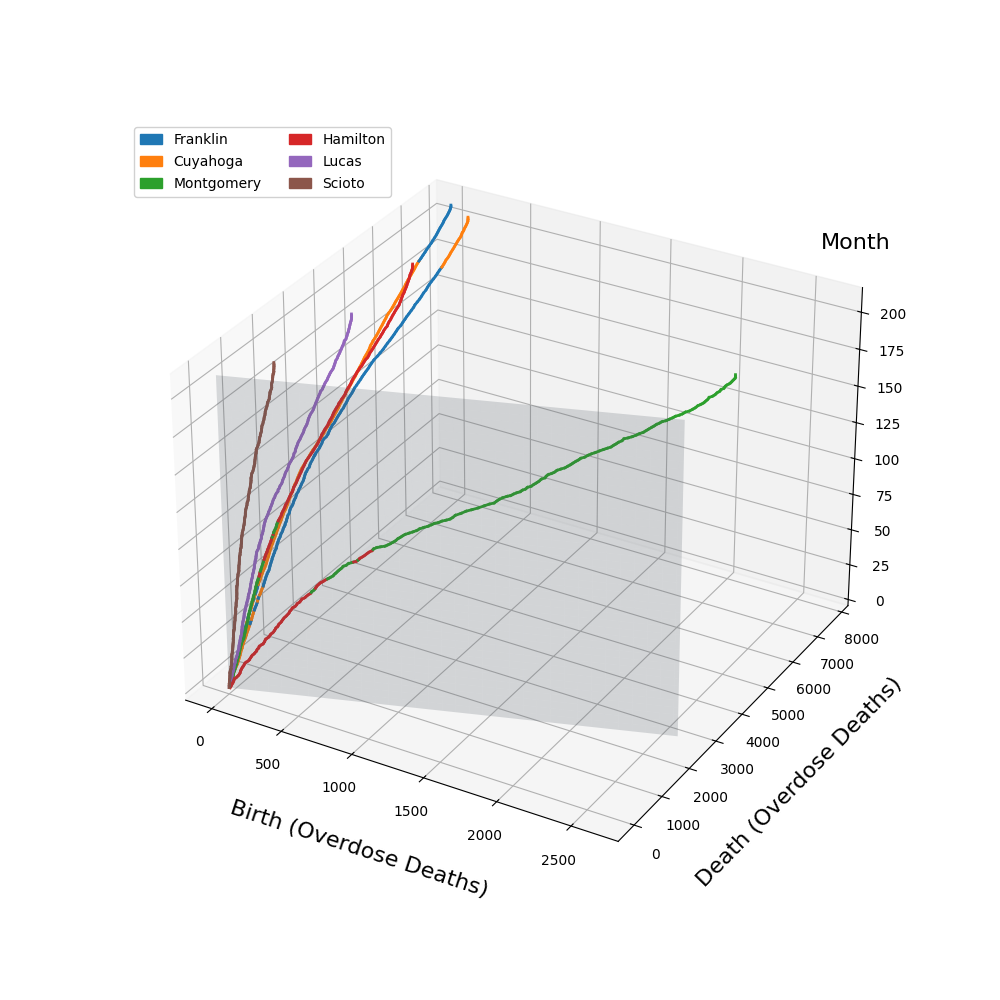}
  \end{minipage}\hspace{0.05\linewidth}
  \begin{minipage}{0.4\linewidth}
    \centering
    \includegraphics[width=\linewidth]{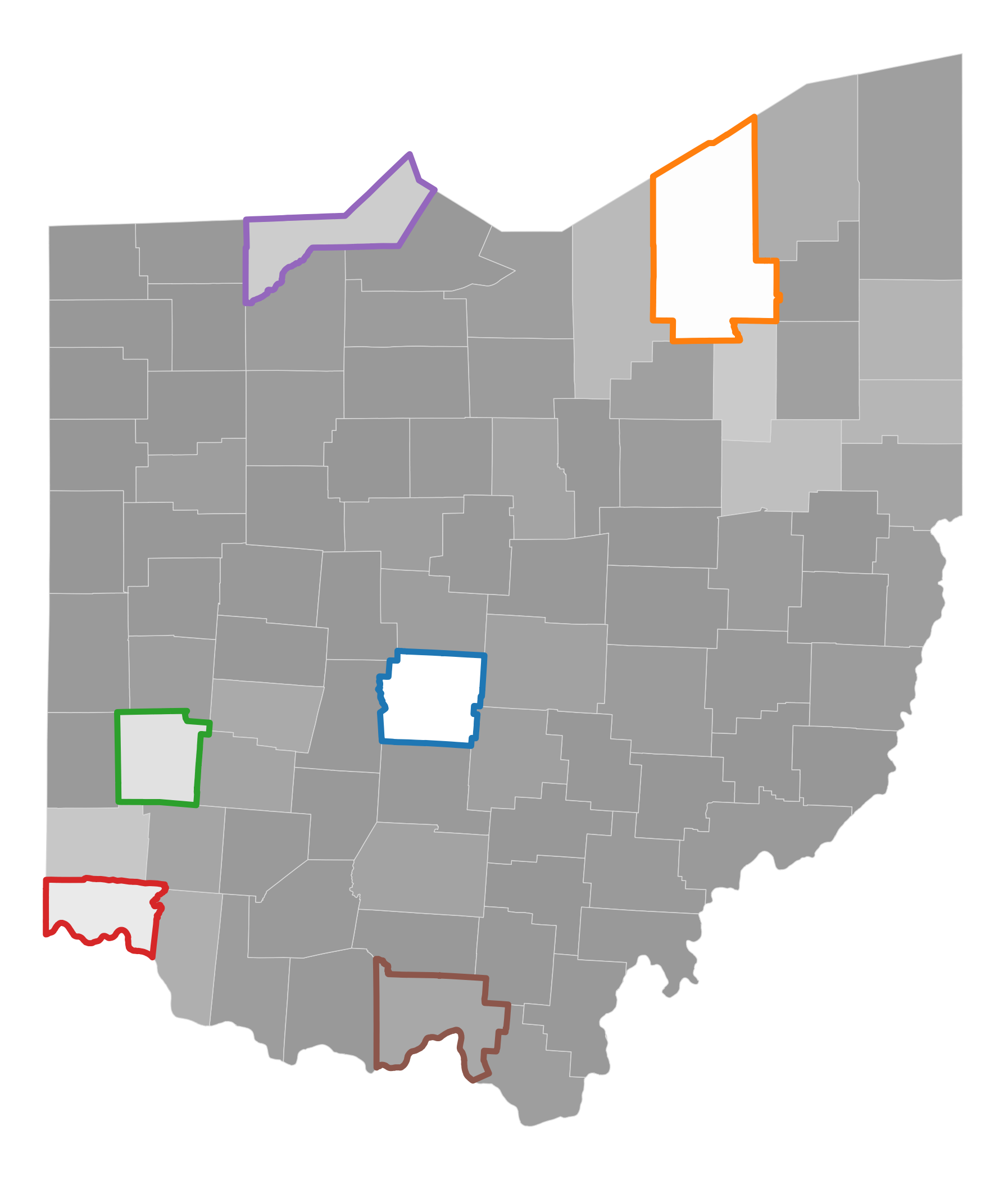}
  \end{minipage}
  \caption{(Left) Vineyard of the six most persistent vines created using cumulative overdose deaths in the counties of Ohio from January 2007 to September 2024. (Right) A grayscale map of the counties of Ohio with intensity determined by the relative cumulative death count in September 2024. The vines are colored according to their representative counties and the border of each of these counties in the map of Ohio is colored accordingly.}
  \label{fig:teaser}
}

\abstract{
Understanding how spatial patterns evolve over time is a complex task that often arises in the analysis of public health data.  In this work, we investigate the use of vineyards from topological data analysis (TDA) in this setting by applying them to time series data related to the overdose epidemic in the state of Ohio.  We begin by proposing statistical tests that can be used in order to evaluate whether vineyards are a reasonable technique to study a spatiotemporal dataset.  We then apply these tests to the data of drug overdose deaths in Ohio and, finding the data suitable, perform a subsequent analysis using vineyards to visualize the evolution of local maxima of death rates throughout the Ohio overdose epidemic. We conclude by developing statistical methods to quantify the significance and uncertainty of vineyard features and by exploring how vineyard-derived summaries can be used for forecasting.
}


\begin{document}


\firstsection{Introduction}

\maketitle

Drug overdose remains a leading cause of preventable death in the United States. Between 1999 and 2022, annual overdose deaths rose more than sixfold, claiming over 107,000 American lives in 2022 alone. The impact of this epidemic is highly disparate, with geography and demography heavily dictating how individual communities are affected. Some of the most severely impacted areas, including several counties in Ohio, have experienced multiple waves of the crisis. This indicates deep-seated spatiotemporal patterns that must be understood to improve mitigation efforts. In an era of tightening public health budgets, identifying counties which fare worse than their neighbors can help provide critical insights into vulnerable populations and empower states to efficiently allocate resources.

However, tracking the spread of overdose deaths is challenging; it does not follow a simple trajectory, but is instead shaped by shifting drug supplies, local economic distress, demographic changes, and varying public health infrastructures. Traditional spatial modeling approaches, such as spatial regression and generalized linear mixed models, struggle here because they require rigid assumptions regarding linearity, locality, and parametric form. Furthermore, high-dimensional models with numerous county-specific coefficients are often difficult to interpret, frequently obscuring broader structural shifts. This complex landscape calls for flexible, non-parametric tools, like those from Topological Data Analysis (TDA). 

The authors previously utilized a prominent TDA tool, Mapper, to investigate these dynamics in Ohio \cite{OhioOverdoseMapper}. While that paper demonstrated Mapper's ability to visualize spatiotemporal and demographic trends with minimal data specifications, it highlighted a distinct challenge: Mapper is primarily a qualitative visualization tool. It does not fit naturally within a quantitative or statistical testing framework. Additionally, the resulting Mapper graphs were frequently cluttered, occasionally obscuring vital information whether embedded in 2D or 3D space.

To address these limitations, this paper introduces the use of vineyards to construct alternative spatiotemporal visualizations of local maxima of drug overdose deaths that occurred in Ohio between 2007 and 2024. Vineyards allow us to filter features by their persistence, while accompanying numerical and vector summaries ensure no significant structural shifts are overlooked. Unlike methods that simply identify local maxima, the use of a sublevel set filtration records how those maxima are separated by the surrounding landscape, allowing us to distinguish isolated peaks from those that are part of a broader spatial structure and to track these relationships through time. We then bridge the gap between TDA and quantitative analysis by proposing a novel framework of statistical tests. These tests serve a dual purpose: they determine whether a given spatiotemporal dataset is appropriate for vineyard analysis in the first place, and they establish the statistical significance of the features within the resulting vineyard diagrams.

Ultimately, by producing interpretable visualizations that highlight the communities most vulnerable to the overdose crisis, this approach offers actionable insights for future public health planning. Similar spatiotemporal analytics in Rhode Island successfully motivated neighborhood-level, geographically targeted interventions \cite{marshall2022preventing}. By identifying Ohio counties that fare worse than their immediate surroundings, we hope to provide a similar blueprint for localized intervention and resource allocation. 

\subsection{Outline and Contributions}

In \hyperref[section:Background]{Section 2} we give some background on vineyards as a tool for studying time series data and highlight some of the complexities in interpreting what each vine represents. We then discuss some of the attempts to develop statistical tests for topological data analysis tools and the contexts in which they have been applied.

In \hyperref[section:StatsforSuitability]{Section 3} we propose our own statistical tests to determine whether vineyards are an appropriate TDA tool to apply to a given dataset. We outline certain null hypotheses that we believe would disqualify vineyards from being an appropriate tool for a dataset, discuss statistical tests to determine whether these null hypotheses hold, and apply those tests to show that the null hypotheses are rejected for our given dataset. In \hyperref[section:StatsforSuitability]{Section 3} we focus on classical statistical tests, and then in \hyperref[section:StatsforSignificance]{Section 6} we introduce TDA-based statistical tests for the same hypotheses. Several of these TDA-based tests, including testing if the scalar field (e.g., drug overdose deaths) is uniform over all counties, and testing for spatiotemporal covariance, appear to be new.  This appears to be the first time that vineyard data has been used as a test statistic for a null hypothesis test.

In \hyperref[section:Methods]{Section 4} we outline the methods that we use to construct our vineyard diagrams from our time series data and in \hyperref[section:Results]{Section 5} we present the results of our vineyard analysis.  We then compare these results to previous work that has been done to analyze the Ohio overdose epidemic.

In \hyperref[section:StatsforSignificance]{Section 6} we develop statistical inference for vineyards, including confidence intervals for numerical summaries, confidence bands for persistence diagrams, and confidence tubes for vines. We use these to create statistical significance tests for vineyard features. We also develop TDA-based hypothesis tests for the null models introduced in \hyperref[section:StatsforSuitability]{Section 3}.

In \hyperref[section:Forecasting]{Section \ref{section:Forecasting}} we explore forecasting of statewide overdose deaths and vineyard-derived summaries of the county-level overdose surface. We forecast the maximum $H_1$ persistence and the county associated with the strongest $H_1$ feature, and discuss how these forecasts could inform local health departments.

Finally, in \hyperref[section:Conclusion]{Section \ref{section:Conclusion}} we summarize our work, discuss its limitations, and propose future research directions to be explored.

The main contributions of this paper are:

\begin{enumerate}
 	\item An application of vineyards to drug overdose data that provides new visualizations of and insights into the dynamics of the Ohio overdose epidemic.
 	
 	\item The development of TDA-based statistical tests to use alongside traditional methods from spatial and spatiotemporal statistics in order to determine the suitability of vineyards for analyzing a dataset.
 	
 	\item A statistical inference framework for vineyards, including confidence intervals, confidence bands, confidence tubes, and null hypothesis tests for the statistical significance of vineyard features.
 	
 	\item Forecasting methods in the context of vineyards.
\end{enumerate}

\subsection{Related Works}

Previous works have investigated the spatial and temporal dynamics of the Ohio overdose epidemic. For example, \cite{bci} created a model for the monthly number of deaths $D_t$ based on its own history. They found statistically significant evidence that $D_t$ depends on the lagged time series $D_{t-1}$, i.e., exhibits temporal autocorrelation. This analysis treated all of Ohio as one block. An alternative, \cite{rosenblum2020rapidly}, applied a generalized linear mixed model (GLMM) to opioid overdose death counts $D_{c,t}$ in each (county, month) pair in Ohio. While this method is general enough to account for both spatial and temporal dependence, since $D_{c,t}$ can be a function of any other $D_{c',t'}$, it faces challenges due to the large number of parameters involved. In \cite{rosenblum2020rapidly}, the potential for one county to affect its geographical neighbors is not explored and linear relationships between each $(c,t)$ and $(c,t-h)$ are assumed.

Several strands of research have focused on the spatial autocorrelation (dependence of $D_{c,t}$ on other $D_{c',t}$) of the drug overdose epidemic in Ohio. Andrew Curtis and members of the Begun Center for Violence Prevention have fit spatial models for drug overdose data in the Cleveland area \cite{curtis2025using, mcmaster2025drug, noriega2023case}, at the census block level. Related work in Cincinnati has employed exploratory spatial statistics, including Local Indicators of Spatial Association (LISA), to identify local overdose clusters \cite{choi2022spatial}.
This kind of technique can produce heat maps and cartographic maps showing which areas are most at risk of overdose spikes, e.g., showing movement of the epidemic into African American neighborhoods by comparing heat maps in one year with the next year. However, the statistical models do not include the time dimension, and we are unaware of how these models can be used for forecasting. Additionally, Adam Eck and his students at Oberlin College use machine learning models (e.g., random forests, gradient boosting, individual decision trees, SVMs, neural networks) to predict county-level overdose death rates \cite{eck}, with the explicit aim of helping guide public policy and resource allocation.

In addition, it is possible to approach spatial and spatiotemporal autocorrelation using a Bayesian framework. Kline, Hepler, and their students have employed Bayesian statistical models to estimate spatial autocorrelation in opioid overdose deaths across Ohio counties, providing insight into geographic clustering and county-level risk factors \cite{hepler2019latent, kline2021estimating, kline2021multivariate}. These papers fit generalized spatial factor models, and look at the relationship between treatments for substance use disorder and drug overdose deaths, in each county. Their algorithm produces spatial weights for each county, interpreted as the degree of unmeasured heterogeneity across counties, causing statistically significant differences that the model cannot explain. This work was extended to add a temporal dimension by Ji \cite{ji2019joint}. Others have fit similar Bayesian spatiotemporal models in the Cincinnati area \cite{li2019suspected}, at the census block level.

Some work uses machine learning or statistical surveillance methods to detect or forecast emerging overdose patterns, including using Gaussian-processes \cite{neill2018machine}, Bayesian logistic growth models to predict future county-level opioid overdose mortality in North and South Carolina \cite{sumetsky2021predicting}, forecasting future opioid-incidence heat maps \cite{choudhuri2019predicting}, monthly opioid-incident forecasting for rapid public-health response in Kentucky \cite{mullen2025forecasting}, point-process models to predict overdose hotspots from heterogeneous EMS and coroner data \cite{liu2021point}, and spatiotemporal neural networks to forecast opioid-overdose from crime data \cite{ertugrul2019castnet}. These papers forecast overdose burden, but they do not forecast persistent-homology summaries of overdose surfaces. There is also work using TDA in forecasting or early-warning settings, including persistent-homology features for Zika forecasting \cite{soliman2020ensemble}, persistence landscapes as early-warning summaries for financial crashes \cite{gidea2018topological}, and persistence-vineyard information to predict future qualitative behavior in tumor-immune simulations \cite{yang2025topological}. Our forecasting analysis is closest in spirit to this second group, but instead of using topology only as an input feature, we forecast a topological summary of the future overdose surface itself.

There has also been work done at the national level and in other states, e.g., \cite{stewart2017geospatial}. The most advanced appears to be Rhode Island, where academic researchers have teamed up with the state health department to develop the PROVIDENT system \cite{marshall2022preventing}. This system uses both machine learning algorithms and statistical models (e.g., spatiotemporal Gaussian processes) to predict overdose mortality at the census block level using SUDORS data (explained in \cite{sudors-ohio}).
The state health department uses these predictions to optimize their deployment of overdose prevention resources at the neighborhood level.

Beyond statistical models, it is also possible to model spatiotemporal spread using Hawkes processes and other methods from dynamical systems.
The middle author used these models to determine the spatiotemporal spread of protests in the USA \cite{rodriguez2023analysis} and in Ukraine \cite{bahid2024statistical}.

Numerous previous papers have applied TDA to other epidemics including the spatiotemporal spread of Covid-19 \cite{PorterVineyards,chen2021topological, ault2022comparison}, Zika \cite{lo2018modeling, soliman2020ensemble, rudkin2023spatial}, influenza \cite{costatopological}, and other contagious diseases \cite{taylor2015topological}. The only investigation into the Ohio overdose epidemic using TDA tools that the authors are aware of is their own, \cite{OhioOverdoseMapper}, which adapted and extended the methodologies from \cite{chen2021topological} to create informative visualizations that identified time delayed correlations between demographic features of Ohio counties and spikes in death outcomes and identified certain communities that were most disparately affected.  

In this paper, we have followed some of the methodologies from \cite{PorterVineyards} to create clearer visualizations of local maxima of overdose death rates and have adapted the work of several authors to create confidence intervals and null hypothesis tests for our context. Our statistical methodology builds on several strands of work in TDA inference. Stability of persistence diagrams implies that small perturbations of the filtering function produce small perturbations of the resulting persistence diagram in bottleneck distance \cite{cohensteiner2007stability}. This stability theorem underlies confidence sets for persistence diagrams, including the bootstrap methods developed in \cite{fasy,chazal_bootstrap}. A related approach converts persistence diagrams into functional or vector summaries, such as persistence landscapes, so that classical statistical procedures can be applied \cite{bubenik,chazal_landscapes}. Recent work has also used persistent homology to test for spatial dependence \cite{SamuelByers2023}. Our setting differs from these applications because the Ohio county complex is fixed and uncertainty enters through a noisy scalar field of county-month overdose deaths, rather than through repeated point-cloud samples. We therefore adapt these ideas to produce parametric simulation-based null hypothesis tests, confidence bands, and confidence tubes for vineyard summaries, based on a fitted negative-binomial count model \cite{hilbe2011negative}.

\section{Background}
\label{section:Background}

\subsection{Vineyards}
\label{subsection:Vineyards}

Persistent homology is one of the foundational tools of TDA and has been modified by various authors to create extended persistence \cite{Cohen-Steiner2009}, zigzag persistence \cite{Carlsson2010} and used on spatial datasets through the persistent homology transform \cite{PHT} and extended persistent homology transform \cite{Turner2024, XPHTSD}. 
The variant of Persistent Homology we will focus on in this paper is that of vineyards, introduced by Cohen-Steiner, Edelsbrunner and Morozov in \cite{VinesandVineyards} where they used vineyards to analyze the dynamics of protein folding.

The idea of vineyards comes from stability results for persistence diagrams which imply that a homotopy of ``tame'' functions, $f_t:X\to \mathbb{R}$, where $X$ is a topological space, induces a continuous path in the space of the persistence diagrams arising from the sublevel set filtrations of $f_t$. We call this path a vineyard.  The tameness condition we are required to satisfy is that each function $f_t$ has only finitely many homological critical values, which simply means that the Betti numbers of the sublevel sets induced by $f_t$ only change finitely many times.  In our case we are always working with functions on finite simplicial complexes which clearly satisfy this condition. In order to apply vineyards, which require continuous information, to real world data where we are often working with a discrete time series, $f_{t_i}$ for $i=1,...,n$, we simply perform a linear interpolation between the functions at each time step and look at the induced vineyards as is done in \cite{VinesandVineyards, PorterVineyards}.

Other researchers have created ad hoc methods for computing vineyards by simply linking points in subsequent persistent diagrams that were sufficiently close together \cite{Ciocanel2021}. This is computationally efficient but comes with a potential loss of interpretability of the vines as they are no longer directly related to the functions at each time step.  As interpretability of the visualizations produced is key to our work in this paper we follow the methodology of \cite{VinesandVineyards, PorterVineyards}.

Methods of analyzing vineyard diagrams usually focus on isolating the most persistent vines \cite{PorterVineyards, Ciocanel2021} because if all vines are plotted the diagrams can become cluttered and hard to interpret. To identify the most persistent vines, we follow the methodology of \cite{PorterVineyards} and rank the persistence of vines by averaging their distance from the diagonal.  
Once we have identified the most persistent vines, we are left with the challenge of interpreting the topological information. In the context of regular persistence diagrams, this is usually done by tracking the birth and death simplex associated to each point in the persistence diagram. This allows us to connect the persistence of a topological feature to location features of our simplicial complex and the associated function values on the birth and death simplices.  Unfortunately, this does not cleanly extend to vineyards due to the way the interval decomposition of our persistence diagram changes over time. Technical details of this problem concerning vineyard modules and algorithmic ways to change the bases over time are discussed in \cite{turner2023representingvineyardmodules} but the basic idea can be illustrated with a simple example.

\begin{figure}
    \centering
    \includegraphics[width=\linewidth]{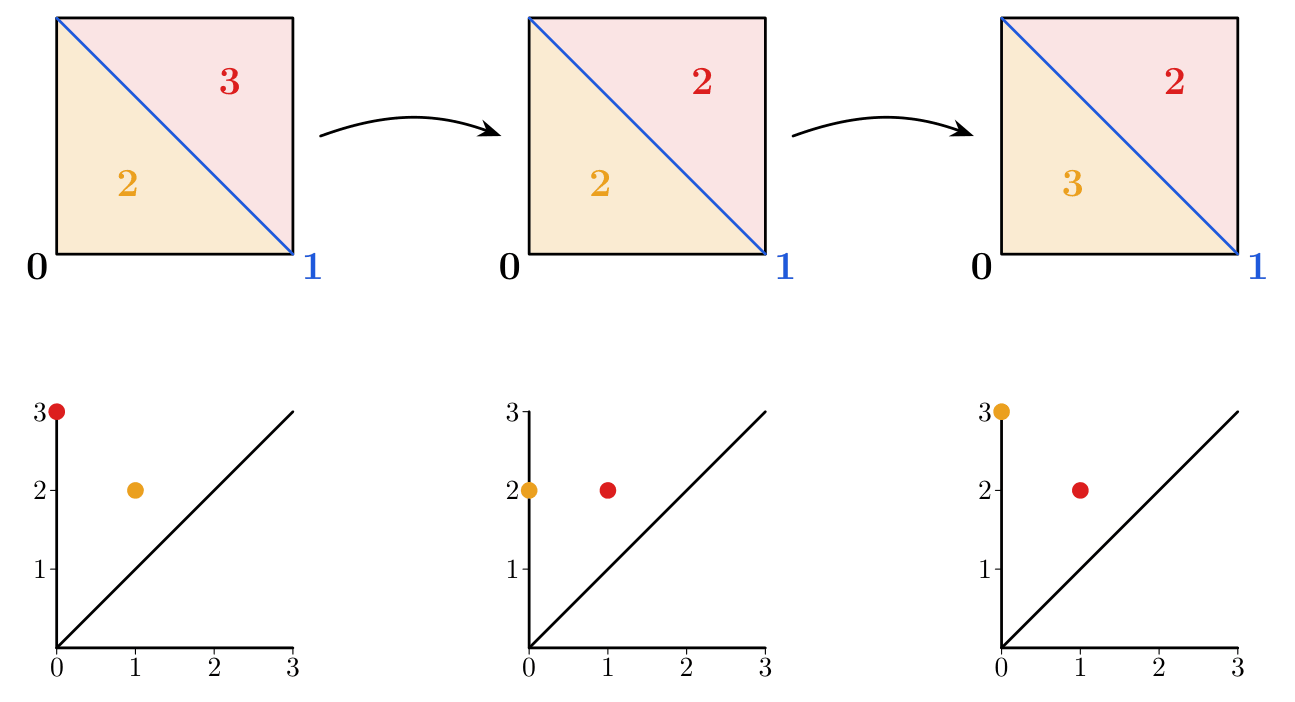}
    \caption{(Top) A sequence of three different integer functions on a simplicial approximation of a square. (Bottom) Corresponding $H_1$ persistence diagrams created by filtering each simplex according to these functions. The persistence classes in each diagram are colored so that they match the color of their death simplex in the square above.  }
    \label{fig:vineswap}
\end{figure}

In \cref{fig:vineswap} we see a simplicial approximation of a square and three snapshots of a homotopy of functions on the simplicial complex.  Beneath these three snapshots of time are the corresponding $H_1$ persistence diagrams for the filtrations defined by the functions at each time step. These diagrams are colored according to the color of their death simplex in the picture above.  

\begin{figure}
	\centering
	\includegraphics[width=0.6\linewidth]{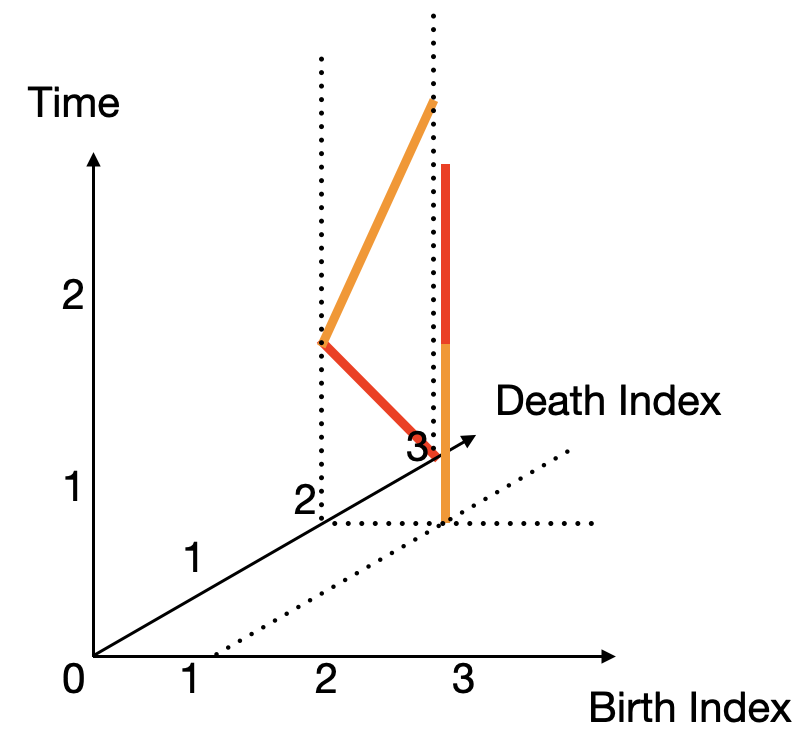}
	\caption{The vineyard of $H_1$ Persistence diagrams produced by linearly interpolating between the snapshot functions of \cref{fig:vineswap}.  Regions of vines are colored by death simplex associated to the persistence point at each point in time.}
	\label{fig:ExampleVineyard}
\end{figure}

Linear interpolation of the functions at these timesteps will induce a vineyard with two vines, one stationary at the coordinate $(1,2)$ over time and the other which goes up and down the plane birth time = 0 as displayed in \cref{fig:ExampleVineyard} where again the vines are colored by the death simplex of the vine at each point in time.

What we observe is that the death simplex changes over time without the vines or the simplices involved needing to be near each other for the swap to occur.  All that is required in this case is that the death values of the two vines are the same at some time in our constructed homotopy.  

In order to interpret the vineyard diagrams we create, we must label them by the birth/death simplices that occur at each time and determine what changes in our setup to cause changes in the birth/death simplices of each vine. This analysis will always be context dependent when studying general vineyards and we discuss what they mean in our context in \hyperref[subsec:Interpreting_Vineyards]{Section 5.2}

The algorithm introduced in \cite{VinesandVineyards} takes worst case linear time in the number of transpositions of the ordering of the simplices. The number of required transpositions is then $O(n)$ where $n$ is the number of simplices. For our setting the number of simplices is a scalar multiple of the number of counties, which is sufficiently small to not create run time issues.

\subsection{Statistics for Persistent Homology}
\label{subsection:StatsforPH}

As observed in \cite{Ciocanel2021}, one of the key unresolved challenges in TDA is understanding how uncertainty in random data relates to uncertainty in the resulting TDA outputs, e.g., determining whether persistent features are statistically significant, and constructing confidence sets (generalizing confidence intervals) to quantify the expected variability.

Statistical inference in TDA is often formulated for persistence diagrams arising from a random point cloud sampled from an underlying distribution or manifold. Much of this literature studies the resulting probability laws, estimators, null hypothesis tests, and confidence sets. A concise and very useful overview is given in \cite{chazal_michel}. Confidence sets for persistence diagrams were developed in \cite{fasy}, where bootstrap and subsampling methods estimate bottleneck-distance uncertainty between an observed diagram and simulated diagram. Features whose distance to the birth = death line are greater than what can be explained by random variation alone are considered statistically significant. The cutoff can be visualized as a band around the diagonal in the persistence diagram.

A second line of work turns diagrams into vector-valued or function-valued summaries, then applies classical statistical procedures. Persistence landscapes were introduced in \cite{bubenik}, with limit theorems that support averages, confidence intervals, and hypothesis tests for functional summaries of diagrams. Related convergence and bootstrap results for landscapes and silhouettes appear in \cite{chazal_landscapes}. Many tests in this area are, in spirit, TDA analogues of two-sample tests: one compares two collections of diagrams, landscapes, or vectorized summaries and asks whether they appear to come from the same distribution. Examples include tests based on distances between diagrams \cite{robinson_turner} and tests based on vectorized persistence diagrams \cite{moon_lazar}.

There has been substantially less work focused on independence, dynamics, or spatial structure. Block-sampled Monte Carlo tests for persistent homology of time series were developed for simulated fMRI data in \cite{abdallah_salch}, bringing TDA into the world of time series statistics. Similarly, persistent-homology tests for spatial dependence were proposed in \cite{SamuelByers2023}. 

Our setting differs from these point-cloud and network applications: our data define a noisy scalar field on a fixed geographic complex, rather than a point cloud. We therefore develop new statistical tests for a variety of purposes, including testing for spatial uniformity, spatiotemporal autocorrelation, and statistical significance of vines. As far as we know, this is the first paper doing statistical analyses in the context of vineyards. In addition to null hypothesis tests, we produce bootstrap confidence bands for vineyards, confidence tubes for individual vines, and scalar confidence around observed features.

\subsection{Ohio Overdose Data}
\label{subsec:DataDescription}

Our dataset of drug-induced deaths in Ohio comes from DataOhio, which tracks death records from the Ohio Department of Health's (ODH) Bureau of Vital Statistics \cite{dataohio}.
This data was reported monthly and our analysis covers the period from January 2007 to September 2024. As such, our data is a $213\times 88$ dataframe containing monthly deaths in each of the 88 counties of Ohio. It is relevant to note that the actual number of overdose deaths might not match the number in our dataset, e.g., because drug overdose is sometimes unreported as a cause of death on death certificates \cite{buchanich2018effect}. The ODH and the Centers for Disease Control try to correct for this, but missing data remains a potential concern. Additionally, as of when we downloaded the data, all of the 2024 numbers were still considered preliminary counts, subject to change.
We use the yearly county level population count from the Census 10-year estimates \cite{census-pop}. 
Overdose death rates are calculated using these monthly death counts and yearly populations.

\section{Statistical Tests for the Presence of Spatiotemporal Effects}
\label{section:StatsforSuitability}

Before utilizing any data analysis tool it is important to address whether the tool is appropriate for the data. As vineyards are designed for data changing in time and space, we first test our Ohio data for evidence of such change. In this section we give a high-level overview of several statistical tests and their results on our data. In Section \ref{section:StatsforSignificance}, we also apply these tests to the residuals of a model that we fit, showing that those residuals do not exhibit spatial or temporal autocorrelation and hence the model fits well. 

For full details on the tests, the reader is referred to the GitHub mentioned in \hyperref[sec:supplemental_materials]{Supplemental Materials}. We focus in this section on classical statistical tests. After constructing the Ohio simplicial complex and its filtration in Section~\ref{section:Methods}, we return in Section~\ref{section:StatsforSignificance} to TDA-based tests of the same null hypotheses.

We focus first on spatial autocorrelation. If a scalar field on a fixed complex is approximately uniformly distributed, then any apparent spatial peaks are likely to be artifacts of random variation rather than meaningful spatial features.

We first test this uniformity hypothesis on the number of drug overdoses. For every month $t$, we let $D_t$ denote the number of drug overdoses in Ohio in that month. If the null hypothesis were true, then we would expect the number of overdoses in each of the 88 counties to be approximately $D_t / 88$. We can therefore simulate from such a model and see how extreme our data looks compared to the simulated distribution. Below, we explain the test statistics we use. If a test statistic has a known distribution, like the chi-square distribution, we can look up p-values using classical tables. Otherwise, we can calculate empirical p-values as the fraction of the simulations where the test statistic was as extreme as what we observed in the data. Of course, the population in Ohio is not uniformly distributed, and counties with a larger population tend to have more overdose deaths, so we also test the analogous population-weighted null hypothesis. To test this hypothesis, we let $P_{c,t}$ be the population in county $c$ and month $t$, $p_{c,t} = P_{c,t}/\sum_{j} P_{j,t}$ be the proportion of the total population that lives in county $c$, and $R_{c,t} = D_{c,t}/P_{c,t}$ be the overdose death rate.

After both of the above hypotheses have been rejected, we know that the data has peaks, but we do not know if the geography matters. Our third null hypothesis is that the data has no spatial autocorrelation, i.e., the geographic assignment of county death rates is exchangeable. To test this, we simulate hypothetical worlds by randomly assigning the death numbers to different counties, breaking the spatial structure. 
We now state the three hypotheses:

\begin{enumerate}
	\item The deaths in each month are uniformly distributed across each county, i.e.,  $$D_{\cdot,t} \sim \text{Multinomial} (D_t; \frac{1}{88},\ldots, \frac{1}{88})$$
	
	\item The deaths in each month are uniformly distributed across population so total deaths in a county becomes proportional to county population, i.e.,  $$D_{\cdot,t} \sim \text{Multinomial} (D_t; p_{1,t},\ldots, p_{88,t})$$
	
	\item The observed county death rates in each month are spatially exchangeable, i.e., conditional on the multiset
$\{R_{1,m},\ldots,R_{88,m}\}$,
each assignment of these rates to the 88 Ohio counties is equally likely.
\end{enumerate}

We now briefly describe our tests of these hypotheses and the results. For full details, see the \hyperref[sec:supplemental_materials]{Supplemental Materials}.
We first tested hypotheses 1 and 2 using traditional statistical tests including Pearson's chi-squared statistic, the likelihood-ratio goodness-of-fit statistic, total variation (half the $L^1$ distance between observed and null hypothesis county shares), largest county-level standardized deviation, and a weighted variance of relative county rates \cite{agresti2013}. For each statistic we computed a global test by aggregating evidence over all months. These global tests rejected both null hypotheses.

We also computed monthwise post-hoc versions of the same tests. To handle the multiple testing problem, we calculated both empirical p-values and Benjamini-Hochberg adjusted q-values within each null-hypothesis/statistic family \cite{benjamini_hochberg_1995}. After this adjustment, the equal-county null was rejected in all 213 months, while the population-weighted per-capita null was rejected in most months.


To test null hypothesis 3, we used Moran's I statistic \cite{anselin2023}. The test rejected hypothesis 3 globally, and rejected a large number of monthly post-hoc tests. 

All of the hypotheses above focus on spatial autocorrelation but not temporal. In previous work, the middle author showed that drug overdose data has temporal autocorrelation \cite{bci}. We now test whether our data exhibit nonseparable spatiotemporal covariance. The null hypothesis is that the spatiotemporal covariance decomposes into a product of spatial and temporal parts, as it would if the two parts were independent. We now state the null hypothesis, letting $Cov(X,Y)$ denote the covariance of the random variables $X,Y$.

\begin{enumerate}[resume]
	\item The spatiotemporal covariance is separable, i.e., $Cov(D_{c,t}, D_{c',t'}) = C_{s}(c,c')C_{t}(t,t')$ for some functions $C_s, C_t$.
\end{enumerate}

As above, we can test this on both raw death counts and deaths per capita. The idea is to simulate a large number of datasets satisfying the null hypothesis, calculate test statistics, and report an empirical p-value.  

We first converted the county-month data into a space-time lag covariance table. County pairs were grouped by graph distance in the county adjacency graph, using bins corresponding to the same county, adjacent counties, graph distance two, and graph distance three or greater. Month pairs were grouped by temporal lag. For each spatial lag $h$ and temporal lag $u$, we computed an empirical covariance $C_{\mathrm{emp}}(h,u)$. Under the separability null hypothesis, this binned covariance surface should be well approximated by a separable product. We therefore fit the best rank-one separable approximation $C_{\mathrm{sep}}(h,u)=\widehat a_h\widehat b_u$ and used the relative residual $\frac{\|C_{\mathrm{emp}}-C_{\mathrm{sep}}\|_F}{\|C_{\mathrm{emp}}\|_F}$ as the test statistic. 

We tested several null models. The first was a global time-block permutation, which reorders blocks of consecutive months. This preserves the spatial maps observed in individual months and retains some local temporal structure within blocks, but disrupts the global chronological ordering. The second was a within-month spatial permutation that randomly reassigns county values within each month, preserving the empirical distribution of rates in each month but destroying the geographic arrangement of those rates. The third was a separable Gaussian matrix-normal model. In this model, we estimate one covariance matrix across counties and one covariance matrix across months, then simulate new county-month scalar fields whose space-time covariance has the product form required by the null hypothesis. Our randomization-based test rejected the null hypothesis under all the null models.

\section{Building Vineyards from County Level Data}
\label{section:Methods}

Having determined that spatiotemporal effects are present in the overdose mortality data, we move to using vineyards for an analysis.
To construct our vineyards, we first create a simplicial approximation of the state of Ohio in which each county is represented by a collection of $2$-simplices. This allows us to extend mortality data defined on counties to a simplicial complex and study its topology through a sublevel set filtration of the mortality data. As the filtration parameter increases, simplices representing counties are added in order of increasing mortality, causing low-mortality regions to appear and merge first. Counties whose mortality values are high relative to their neighbors enter later and can temporarily enclose regions already present in the filtration, giving rise to persistent $H_1$ classes. Consequently, long-lived loops in the filtration correspond to spatially localized maxima in the mortality data.

For our vineyard construction, we use $H_1$ persistence as opposed to $H_0$ persistence of a descending filtration to ensure that all persistence classes have finite birth and death values. This is guaranteed by the topology of the underlying space being homotopy equivalent to a disk and is useful for some of the forecasting described in \hyperref[section:Forecasting]{Section 7}. This methodology is generally applicable to circumstances where the underlying topology is different than a disk, but considerations would have to be made regarding interpretation of essential homology classes of the domain which would depend on the particular application. These considerations would also likely require small modifications to the implementation of statistical tests described in \hyperref[section:StatsforSignificance]{Section 6} that are context specific.

To construct this simplicial approximation, we follow the methodology of \cite{PorterVineyards}, using adjacency information obtained from a shapefile. In this method, each county is first replaced by a polygon with sufficiently many sides to realize all of its adjacencies. These polygons are then glued together along edges corresponding to county adjacencies and triangulated to produce a simplicial complex. Throughout this process, we record which triangles belong to the polygon associated with each county.

Given a function $f$ on the counties of Ohio, we define a function $\hat{f}$ on the simplicial approximation. For each $2$-simplex $\sigma$, we set $\hat{f}(\sigma)$ equal to the value of $f$ on the county whose polygon contains $\sigma$. To define $\hat{f}$ on the $1$-simplices, we distinguish two cases. If $\tau$ is a boundary edge of the simplicial disk representing Ohio, we define $\hat{f}(\tau)$ to be the minimum value of $f$. Otherwise, we set
\[
\hat{f}(\tau)=\min_{{\sigma \mid \tau < \sigma}} \hat{f}(\sigma),
\]
where $<$ denotes the face relation. The values on vertices are defined analogously.

Assigning values in this manner ensures that $\hat{f}$ is a valid filtration function: every simplex enters the filtration no later than any simplex containing it. Consequently, all simplices associated with a county enter the sublevel set filtration simultaneously, so the filtration of $\hat{f}$ faithfully reflects the ordering of counties induced by $f$. Moreover, the resulting simplicial complex retains the topology of Ohio.

Using this framework, we study two vineyards: one using a time series of cumulative deaths in each county of Ohio and one using population normalized cumulative deaths (death rates). The choice of using cumulative deaths rather than raw death counts was to better satisfy the continuity assumption underpinning our usage of vineyards. We additionally explored a version with noncumulative deaths averaged over a window of time for smoothing purposes, as done in \cite{PorterVineyards}, however these approaches communicated similar findings to the cumulative pictures while exhibiting less clarity.  We have included these figures as supplementary materials in the GitHub, see \hyperref[sec:supplemental_materials]{Supplemental Materials}.

To turn our time series data into a continuous collection of functions, which is required for constructing vineyards, we follow the same approach as \cite{VinesandVineyards,PorterVineyards} and perform a linear interpolation between the functions at each time step. We also follow the approach of \cite{PorterVineyards} in coloring the vines according to the county associated to the death simplex of the vine at each point in time. For our visualizations, we only plot the six most persistent vines as measured by average distance from the diagonal in order to reduce clutter and highlight the most important features.

\section{Results}
\label{section:Results}

\subsection{Vineyard Diagrams}

\begin{figure}
    \centering
    \includegraphics[
        width=\linewidth,
        clip
    ]{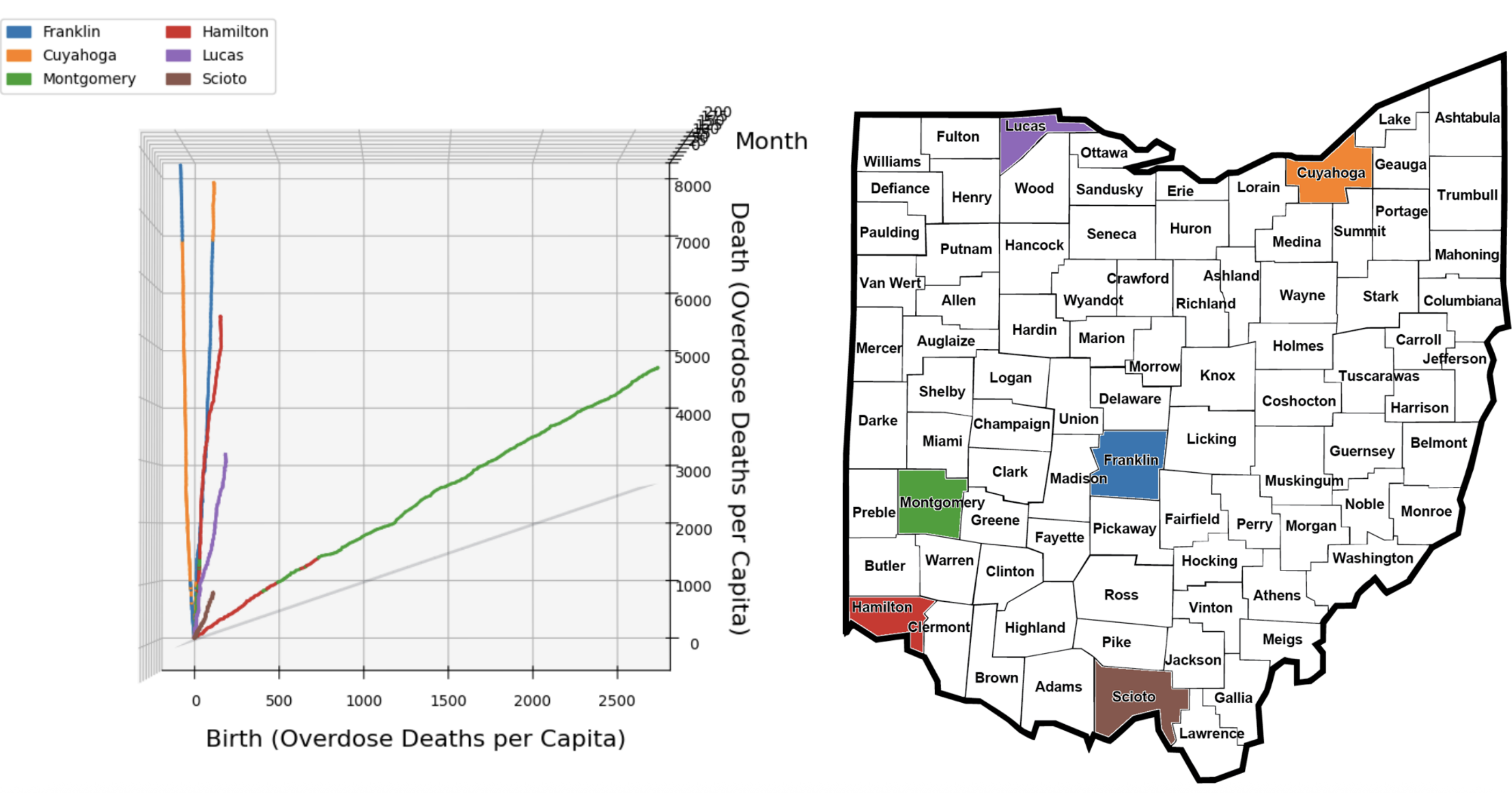}
    \caption{(Left) A top-down view of a vineyard constructed using cumulative overdose death counts across Ohio counties, with colors indicating the county associated with each death simplex at each point in time. (Right) A map of Ohio counties colored according to the legend used in the vineyard diagram.}
    \label{fig:cumulative_vineyard}
\end{figure}

In \cref{fig:cumulative_vineyard} we see the vineyard created by the raw count cumulative overdose data. Adjacent is a map of the counties of Ohio colored according to the legend in the vineyard diagram, which indicates the death simplex of each vine at each time.  A similar diagram is presented in \cref{fig:teaser} with the outline of the counties colored according to their legend in the vineyard diagram, while the counties themselves are colored on a map depicting the relative level of cumulative overdose deaths at the end of the temporal window studied. A .gif file showing the temporal evolution of the cumulative overdose-death map is available on our GitHub, see \hyperref[sec:supplemental_materials]{Supplemental Materials}.


Here we can see that the most persistent vines are typically those that contain major cities in Ohio, notably Franklin county containing Columbus, Cuyahoga County containing Cleveland, Hamilton county containing Cincinnati, Montgomery county containing Dayton and Lucas county containing Toledo. The sixth most persistent vine in this diagram, which stays relatively close to the diagonal, is Scioto county.
Scioto contains the town of Portsmouth, which is known to be one of the hardest hit rural areas of Ohio \cite{pillmill}.  This is very similar to the findings of \cite{OhioOverdoseMapper} which found Franklin, Cuyahoga, Hamilton, Lucas and Montgomery County to be important when considering cumulative overdose deaths, but only found Scioto to be prominent when considering population normalized data. This suggests that vineyard visualizations are more robust to clutter and differences in scale than Mapper visualizations.

\begin{figure}
    \centering
    \includegraphics[width=\linewidth]{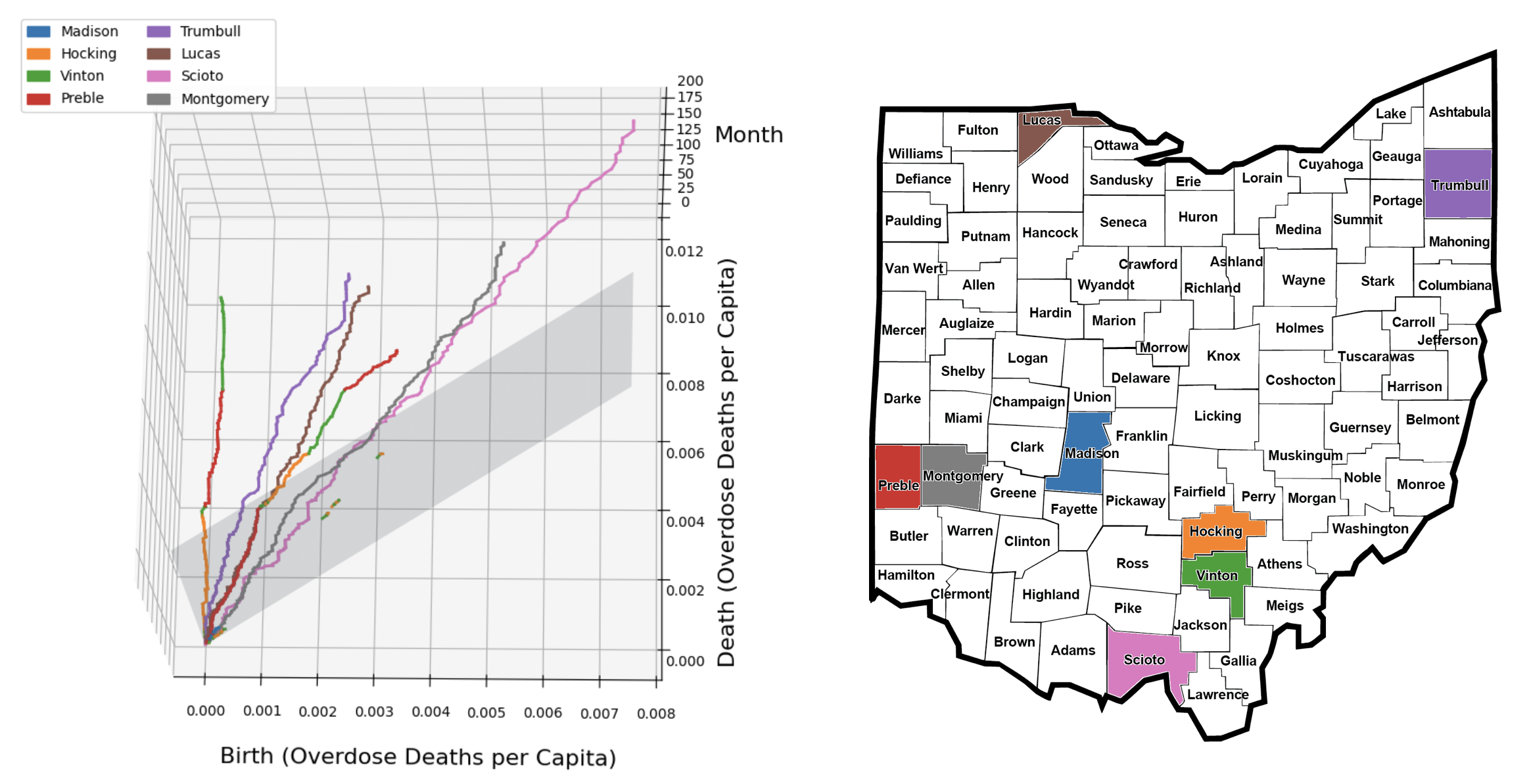}
    \caption{(Left) A top-down view of a vineyard constructed using population-normalized cumulative overdose death counts across Ohio counties, with colors indicating the county associated with each death simplex at each point in time. (Right) A map of Ohio counties colored according to the legend used in the vineyard diagram.}
    \label{fig:normalized_vineyard}
\end{figure}

In \cref{fig:normalized_vineyard} we see the vineyard created by our population normalized cumulative overdose data next to a map of Ohio with counties colored according to the legend in the vineyard diagram. Here we can see that Scioto, Montgomery, Lucas and Trumbull counties contain the death simplex of one vine each that persists throughout the whole time series, and the other two vines contain death simplices belonging to Madison, Hocking, Vinton and Preble County.
These findings go beyond what was discovered in \cite{OhioOverdoseMapper}, which did not identify Hocking, Madison, Preble and Vinton as counties of interest when considering population normalized data.  This suggests once again that these visualizations are capable of picking out counties that are missed due to visual clutter in Mapper visualizations.

\subsection{Interpreting Vineyard Information}
\label{subsec:Interpreting_Vineyards}

\begin{figure}[t]
    \centering
    \includegraphics[
        width=.8\linewidth,
        clip
    ]{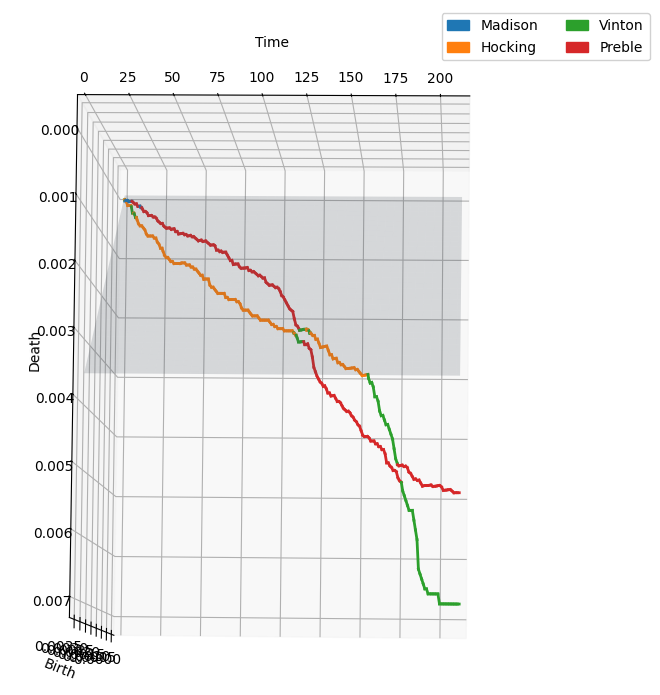}

    \vspace{0.5em}

    \includegraphics[
        width=\linewidth
    ]{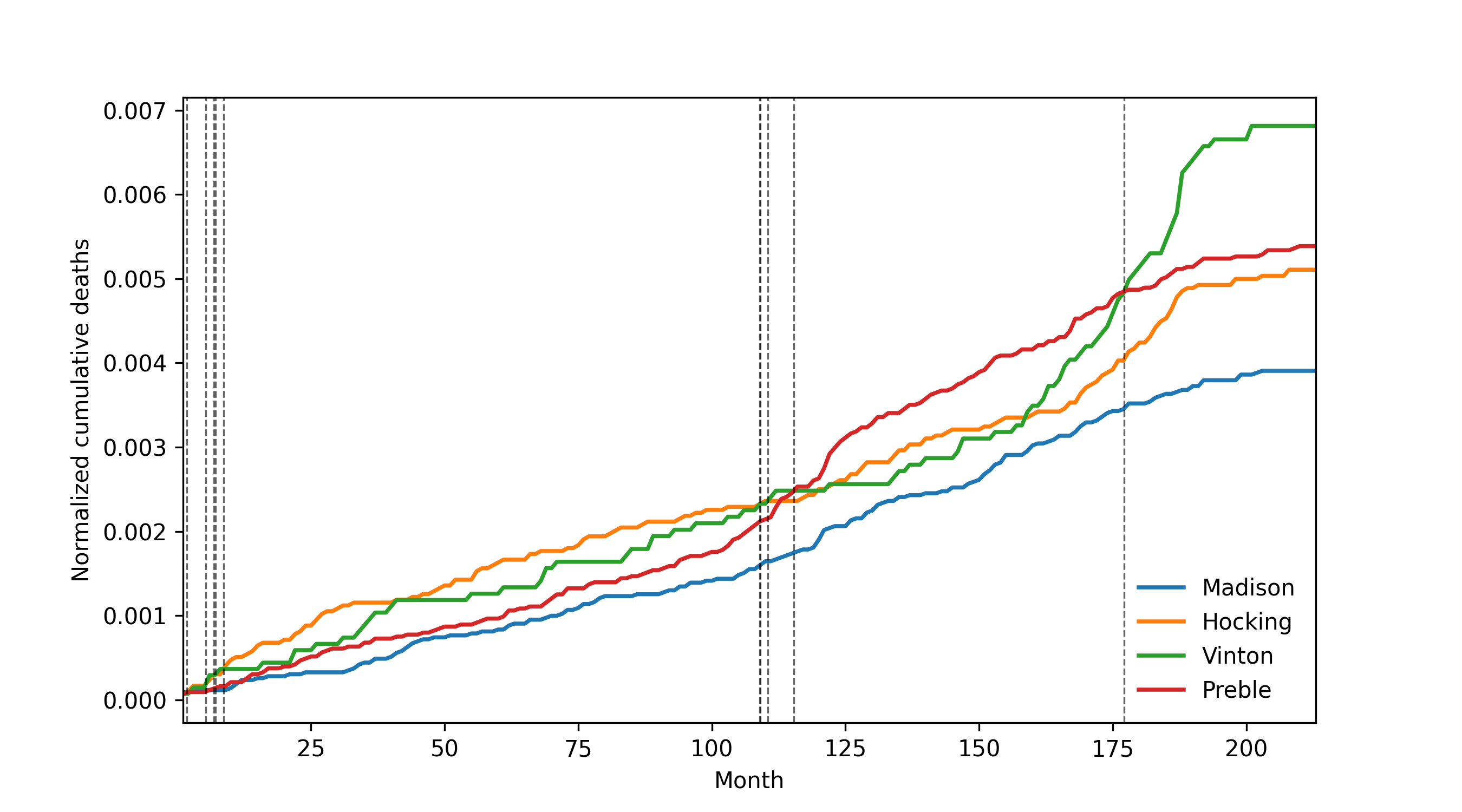}

    \caption{(Above) Vineyard Diagram of the two most persistent vines with death simplices coming from the counties of Madison, Hocking, Vinton or Preble. The vines are colored according to the county associated to the death simplex at each time.  
(Below) Plot of the normalized cumulative deaths of Madison, Hocking, Vinton and Preble county over time.  Dotted lines indicate times that the top two counties swap rank in normalized cumulative deaths, corresponding to a swapping in labels between the vines.}
    \label{fig:hmvp_combined}
\end{figure}

We are now left with how to interpret the information in each vine.  The key subtlety in deriving information directly from the vineyard is that the vines map directly onto a topological feature, but this does not always directly translate to a county or county level information.  We have colored our vines by the death simplex of each persistence class at each time to ensure that the ``death time'' of the vine corresponds directly to the (normalized) cumulative count of overdose deaths in that county at that time.  This is a direct consequence of choosing to focus on $H_1$ persistence of our sublevel set filtrations and tying county information to the 2-simplices that ``kill'' these persistence classes.  This methodology also gives meaning to birth times of our vines as the lowest (normalized) cumulative death count required to form a loop of counties in the sublevel sets that encircle the county identified by the death simplex as a local maxima; i.e. that there is no higher local max encircled by the same counties.  One consequence of this is that the persistence class with the greatest death value has a birth time equal to the minimum over all counties as this is the value that the boundary of Ohio has been assigned in our filtration.  Another is that we can understand times at which death simplices on a vine change, or jump, to a different vine, as being times in which one county has overtaken another as the greatest local max encircled by a set of counties.  This explains why the counties contributing to the death simplices of a vine can be geographically distant from each other when the encircling region defining the birth of the vine is sufficiently large. To explore this idea further, we can analyze  the two vines in \cref{fig:normalized_vineyard} whose death simplices come from Madison, Hocking, Vinton and Preble Counties using the plots in \cref{fig:hmvp_combined}. 

It is clear from analyzing this figure that the times when the death simplex on the most persistent vine swaps with the death simplex on the second most persistent vine coincide with times that the overall normalized death count in one county overtakes another as the highest of this set of counties.  Note that this has nothing to do with the birth times of the two vines which develop based on when counties fill in to create smaller holes in our simplicial approximation of Ohio that isolate our local maxima.  Note also that aside from Hocking and Vinton, these counties are not adjacent geographically and yet represent the same vine. This swapping also occurs for the two vines in \cref{fig:cumulative_vineyard} whose death simplices belong to Franklin and Cuyahoga County, when Franklin overtakes Cuyahoga as having the greatest cumulative death count of all counties of Ohio towards the end of the time window we considered.

Due to the use of cumulative data, the evolution over time tracks the growing difference between the death tolls in a county that determines a local max and the death tolls of their ``neighbors'' as described above. This also means the vines will monotonically increase in both coordinates which explains the shape of the vines in \cref{fig:cumulative_vineyard} and \cref{fig:normalized_vineyard}.

Other insights we can draw from these diagrams are that Scioto appears in both of our vineyard diagrams emphasizing the extent of the impact of the Ohio overdose epidemic in this community and the fact that Madison, Preble and Trumbull are highlighted in \cref{fig:normalized_vineyard} could point to spatiotemporal spread out of major cities to affect regional areas.

\section{Confidence Intervals and TDA-based Statistical Tests}
\label{section:StatsforSignificance}

\begin{figure}[t]
    \centering

    \begin{minipage}[t]{0.49\linewidth}
        \centering
        \includegraphics[width=\linewidth]{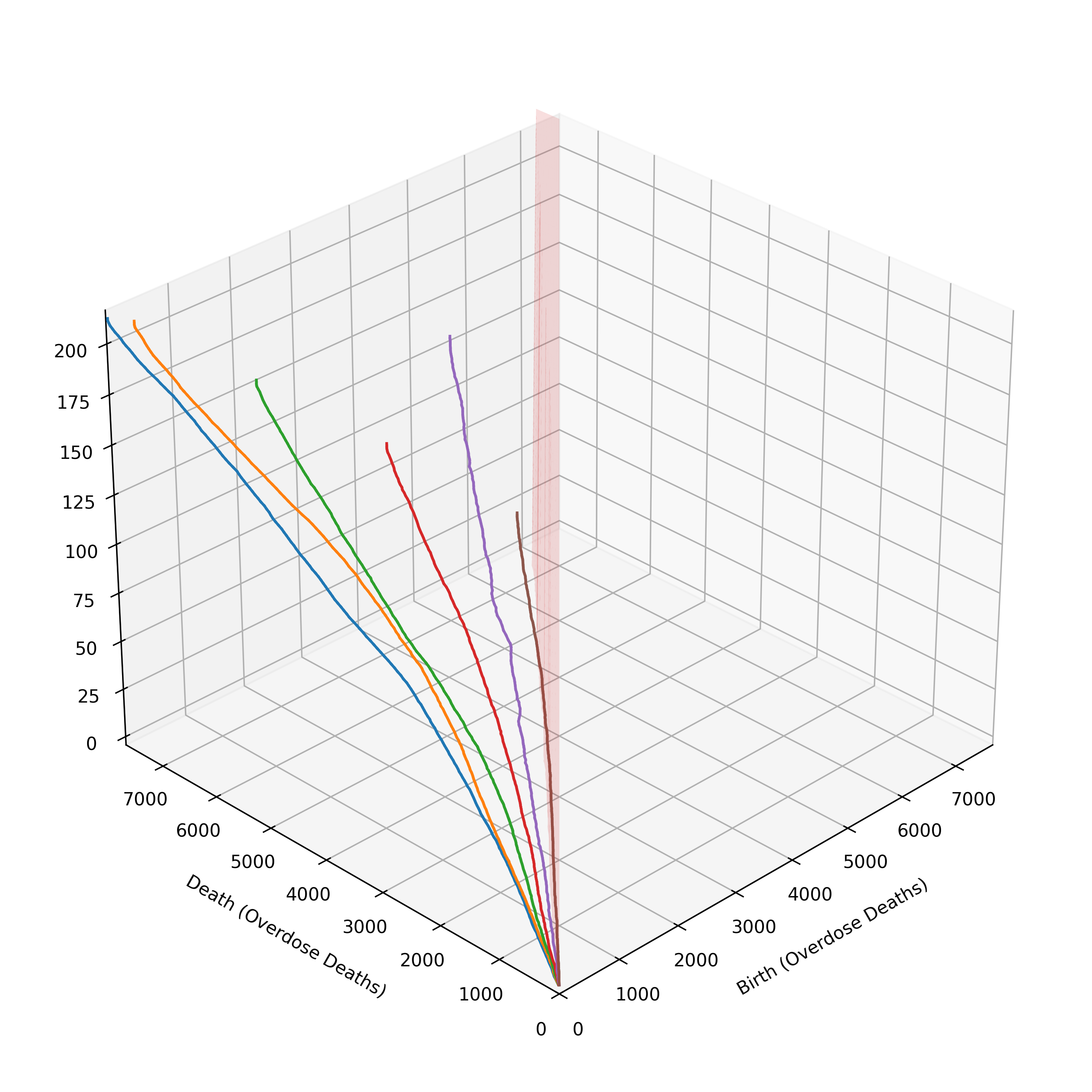}
    \end{minipage}
    \hfill
    \begin{minipage}[t]{0.49\linewidth}
        \centering
        \includegraphics[width=\linewidth]{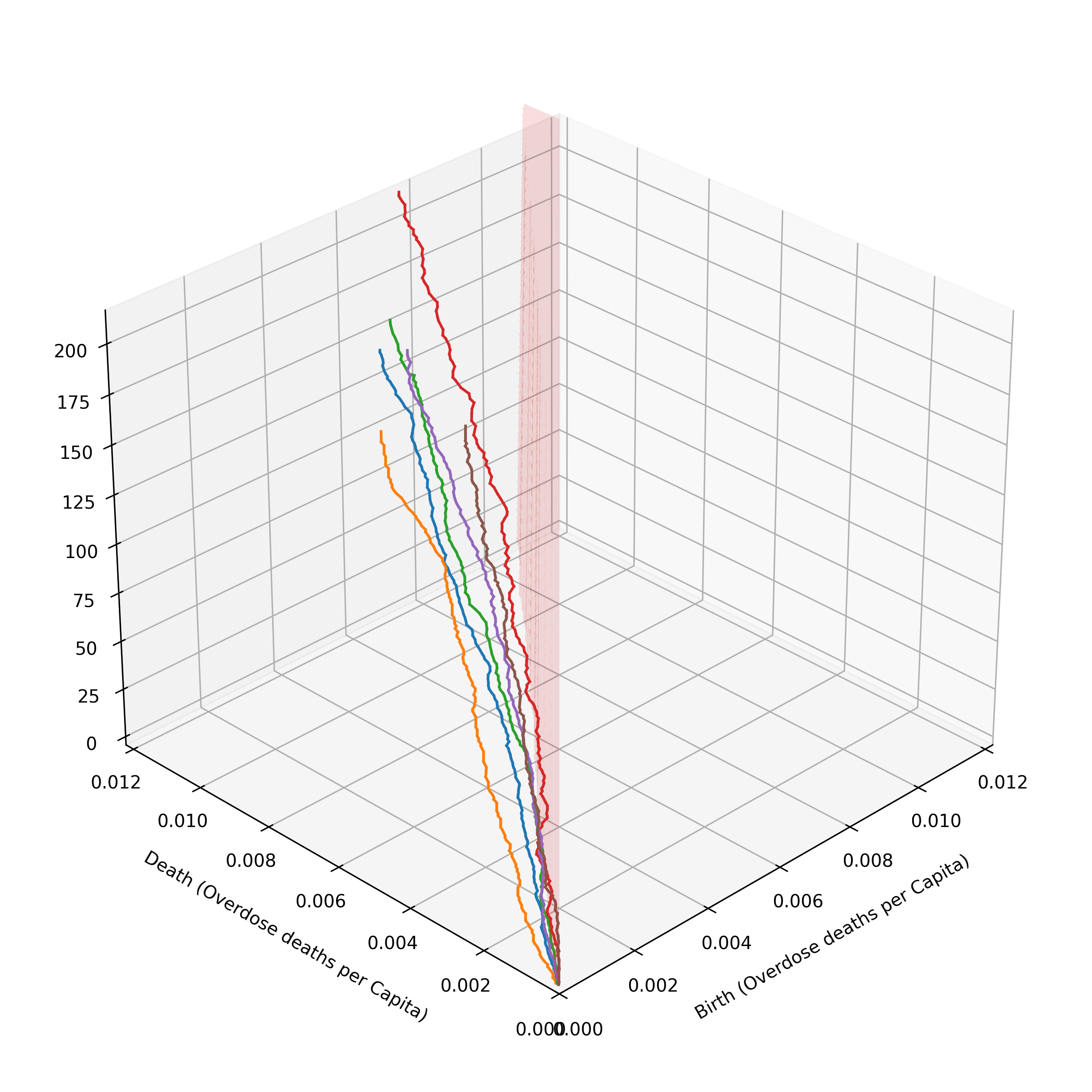}
    \end{minipage}

    \caption{95\% confidence band (shown in red) for the cumulative death count vineyard (Left) and the population-normalized cumulative death count vineyard (Right).}
    \label{fig:confidence_bands}
\end{figure}

In this section, we discuss analogs of confidence intervals in the context of vineyards, tests for the statistical significance of vineyard features, and TDA-based tests for the null hypotheses considered in Section \ref{section:StatsforSuitability}. For full details on the negative-binomial model we used, simulations, code, and resulting confidence intervals, bands, tubes, and statistical tests, please consult the GitHub in  \hyperref[sec:supplemental_materials]{Supplemental Materials}.

\subsection{Confidence intervals, bands, and tubes for vineyards}
\label{subsec:confidence-intervals}

We wish to quantify uncertainty in the vineyard features by bringing the ideas of confidence intervals into the setting of vineyards. To do this, we fit a negative-binomial model \cite{hilbe2011negative} for county-month death counts, including population, county effects, month effects, and county-time interaction structure as covariates. Statistical tests, including those from \hyperref[section:StatsforSuitability]{Section 3}, found no statistically significant spatial or temporal autocorrelation in the residuals, suggesting that the fitted model captures the main spatiotemporal structure. We then generate 1000 bootstrap county-month death tables from this fitted model, viewing these as alternative hypothetical worlds that could plausibly have happened. For each, we convert counts to deaths per capita, recompute the monthly $H_1$ persistence diagrams and vineyard summaries on the same fixed Ohio county complex, and use these bootstrap samples to construct confidence intervals. 

We first construct monthly $95\%$ perturbation thresholds.
This bootstrap distribution produces two complementary uncertainty summaries.
First, we construct monthly diagonal confidence bands for $H_1$ birth-death points following \cite{fasy}. 
The width of each monthly confidence band is determined by the $95\text{th}$ percentile of the bottleneck-distance distribution.
When a persistence point lies
above the band it is considered statistically significant. We adapt this confidence band construction to the context of vineyards by aggregating these bands as they evolve through time to construct a volumetric region, see Figure \ref{fig:confidence_bands}. We say that the vines represent a
statistically significant topological feature during the times at which they extend past this region.

\begin{figure}[t]
    \centering

    \begin{minipage}[t]{0.44\linewidth}
        \centering
        \includegraphics[width=\linewidth]{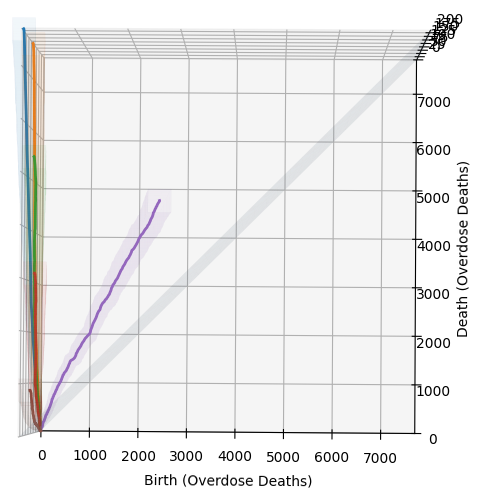}
    \end{minipage}
    \hfill
    \begin{minipage}[t]{0.55\linewidth}
        \centering
        \includegraphics[width=\linewidth]{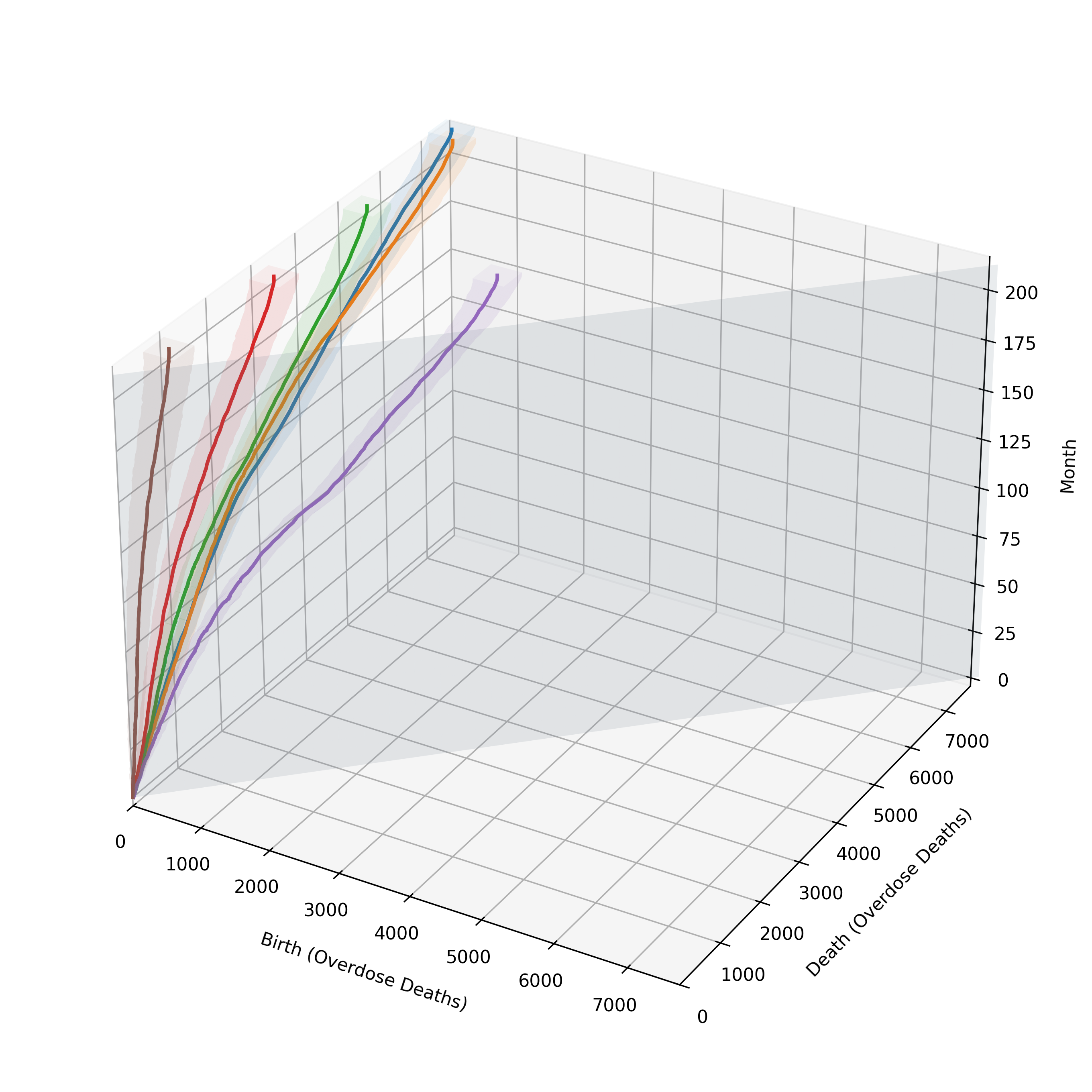}
    \end{minipage}

    \caption{(Left) Top-down view and (Right) angled view of the six most persistent vines created from the cumulative death count data with a 95\% confidence interval tube plotted around them.}
    \label{fig:cumulative_vineyard_Conf_Int}
\end{figure}

\begin{figure}[t]
    \centering

    \begin{minipage}[t]{0.44\linewidth}
        \centering
        \includegraphics[width=\linewidth]{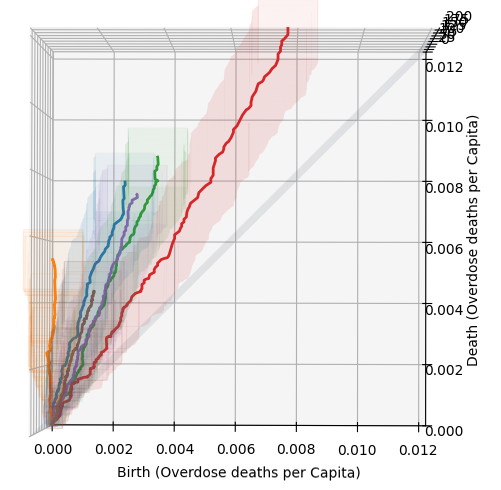}
    \end{minipage}
    \hfill
    \begin{minipage}[t]{0.55\linewidth}
        \centering
        \includegraphics[width=\linewidth]{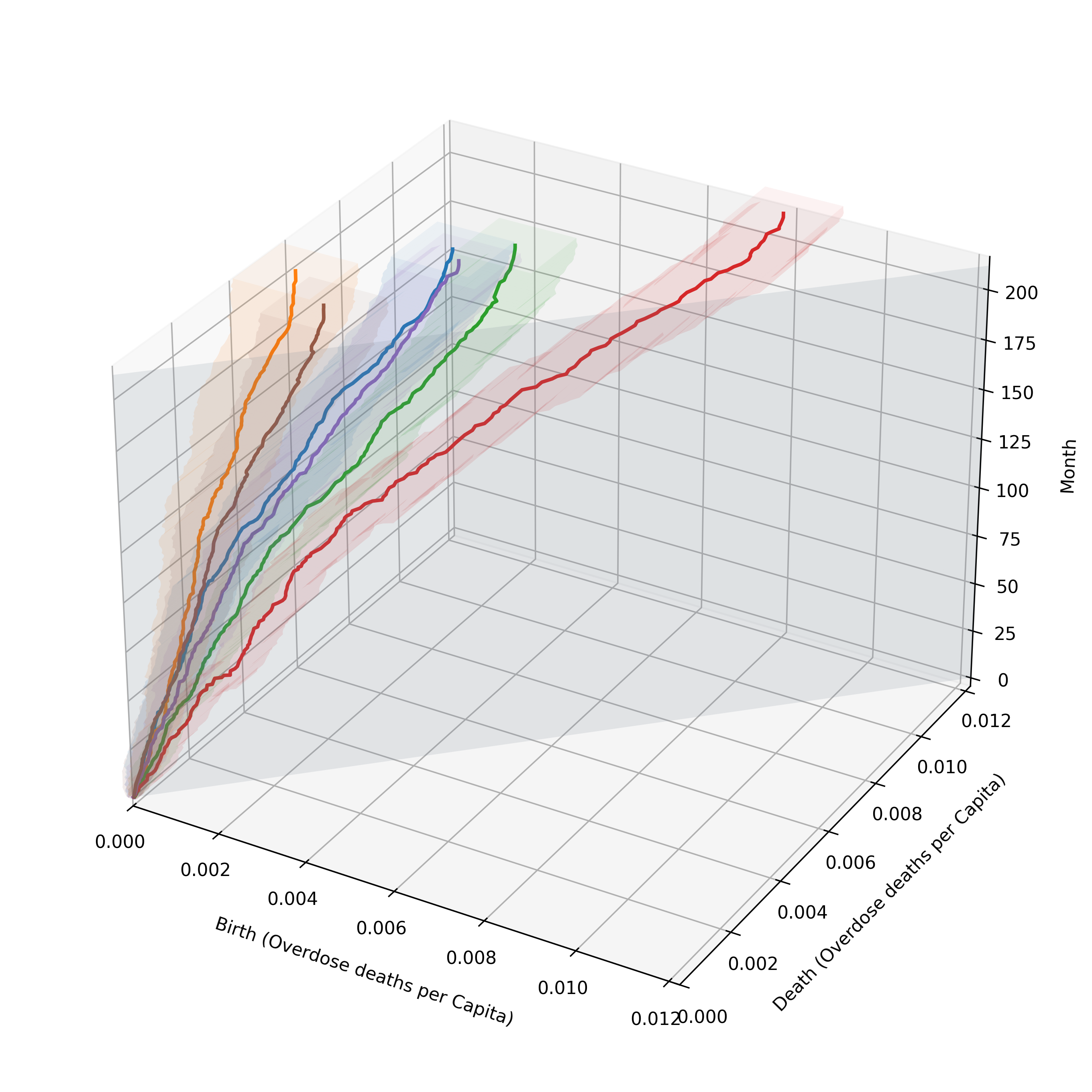}
    \end{minipage}

    \caption{(Left) Top-down view and (Right) angled view of the six most persistent vines created from the population-normalized cumulative death count data with a 95\% confidence interval tube plotted around them.}
    \label{fig:normalized_vineyard_Conf_Int}
\end{figure}

Second, we construct geometric confidence tubes directly around the individual vines shown in Figures \ref{fig:cumulative_vineyard_Conf_Int} and \ref{fig:normalized_vineyard_Conf_Int}. 
Since the bottleneck distance measures the largest $\ell_\infty$ distance between optimally matched diagram points, it naturally establishes an upper bound on the variation of any single persistence feature in the birth-death plane. 
Consequently, the (one-sided) $95\text{th}$ percentile of the bottleneck distance distribution defines a perturbation radius for persistence features in each month. 
By computing these thresholds on a monthly basis and sweeping them along the temporal axis, we generate a continuous tube structure centered on each vine. 
These tubes measure the scale on which $H_1$ features move in the birth-death plane under perturbations of the death counts over time. Vines are statistically significant when their confidence tubes do not intersect the birth = death plane.

In addition to these uncertainty measures, we compute bootstrap confidence intervals for scalar summaries of the diagrams, such as maximum persistence among the strongest $H_1$ features. Finally, we compute confidence intervals for persistence-landscape summaries following \cite{bubenik}. As far as we are aware, this is the first time confidence intervals have been used in the context of vineyards.





\subsection{TDA-based tests for spatial autocorrelation and spatiotemporal effects} \label{subsec:tda-tests}

In \hyperref[section:StatsforSuitability]{Section 3} we tested four null hypotheses using classical statistical methods, to determine if vineyards were an appropriate modeling framework. We now illustrate TDA-based tests for these same null hypotheses.

For hypotheses 1 and 2, we use two $H_1$-based test statistics. For hypothesis 1, the scalar field is the raw county death count. For hypothesis 2, the scalar field is the county death rate per capita. In both cases, we used two $H_1$-based test statistics, computed on both the original data and the data simulated from the null hypotheses. The first records the maximum $H_1$ persistence in each month's persistence diagram. The second is total $H_1$ persistence, i.e., the sum of all persistences in our $H_1$ persistence diagrams for a given month. Here the persistence of a point in our diagram is equal to death time minus birth time. As the distribution of these test statistics is unknown, we calculate empirical p-values. As with the classical tests, we reject both global null hypotheses, and a large number of the monthly post-hoc tests. Indeed, the TDA tests are able to reject hypothesis 2 for certain months that are not flagged as significant by the classical tests.

To test hypothesis 3, we used a separate TDA statistical test motivated by the work in \cite{SamuelByers2023}. The TDA test works by looking at the adjacency graph of counties and analyzing the $H_0$ persistence diagrams created using monthly death data. The choice of $H_0$ persistence over $H_1$ persistence is because spatial autocorrelation is more associated with connected clusters of similar values than loop-like holes. As with hypothesis 2, the TDA test was able to reject certain months missed by the classical test. 

Lastly, for hypothesis 4, we simulated from the same null models as in \hyperref[section:StatsforSuitability]{Section 3}, but now used TDA-based test statistics. We first computed the monthly $H_1$ persistence diagram associated to the sublevel set filtration of the Ohio county complex. We then summarized the resulting time-indexed sequence of persistence diagrams using statistics designed to measure both the magnitude and temporal organization of loop-like spatial features. These included the maximum $H_1$ persistence in each month, the total and average maximum persistence over time, the bottleneck distance between consecutive monthly diagrams, the total bottleneck path length, the mean month-to-month restructuring statistic, and the lag-one autocorrelation of monthly maximum $H_1$ persistence. In addition, for the separable Gaussian TDA test we computed exact vineyards for each simulated dataset and compared vineyard statistics such as the average persistence of the strongest vine and the duration-weighted average persistence of the saved vines. As far as we are aware, this is the first time vineyards have been used as a test statistic for a null hypothesis test, and this is the first TDA-based test for spatiotemporal covariance. In all cases, empirical p-values were statistically significant (even with a multiple-testing correction), so all the TDA tests reject all the null models for the Ohio data.

\section{Forecasting}
\label{section:Forecasting}

One of the most powerful applications of time series analysis is forecasting, i.e., making predictions about what will happen next. In this section, we discuss forecasting of drug overdose data with and without TDA. We find that classical negative-binomial models struggle to forecast because of contemporaneous spatial autocorrelation. Univariate time series are easier to forecast, and we do so for both Ohio-wide deaths $D_t$ and max $H_1$ persistence $M_t$, which simplifies the situation while preserving the key topological information. We then explain how to translate a forecast of $M_t$ back into actionable information for local health departments, to aid them in optimally distributing resources to save lives. As mentioned in  \hyperref[section:Methods]{Section 4}, our 2024 overdose data was preliminary, so we do not include 2024 in either our training or testing data. For full details, see the GitHub in \hyperref[sec:supplemental_materials]{Supplemental Materials}.

The most direct kind of forecasting for public-health applications would be to predict future county-level overdoses: if one could accurately forecast $D_{c,t+1}$ for every county $c$, then local health departments could surge resources to areas where overdoses are predicted to spike. Clearly, the best possible model would involve demographic factors, economic covariates, information about the drug supply, etc. Attempting this would take us too far afield, so we only attempt to predict $D_{c,t+1}$ from other $D_{c',t-h}$ for various $c'$ and $h\geq 0$, i.e., using the past and present. Note that our negative-binomial model from Section \ref{section:StatsforSignificance} cannot be used here because it would need terms of the form $D_{c',t+1}$ to predict $D_{c,t+1}$, i.e., it would still need some future knowledge. 

%
%

We explored county-level negative-binomial forecasting models that predict $D_{c,t+1}$ only from past and present data, but none of these models achieved independent residuals. This is consistent with the statistical tests in Section \ref{section:StatsforSuitability}, which show that the data contain spatial, temporal, and nonseparable spatiotemporal structure that simple county-level forecasting models do not fully capture. 

Instead, we forecast two scalar time series. The first is statewide overdose deaths $D_t$, which is epidemiologically meaningful but ignores county heterogeneity. The second is
\[
M_t=\max_{(b,d)\in PD_t^{(1)}}(d-b),
\]
where $PD_t^{(1)}$ is the finite $H_1$ persistence diagram for the county deaths-per-capita field $R_{-,t}$ in month $t$. Letting $m = \min_c R_{c,t}$ and $M = \max_c R_{c,t}$, we have $M_t = M-m$, because the Ohio complex is topologically equivalent to a disk, the boundary of the disk is assigned filtration value $m$ in the sublevel set filtration, and the filling of that boundary occurs at filtration $M$. For this reason, we can interpret $M_t$ as the range of the county deaths-per-capita field. 

Although $M_t$ admits this simpler interpretation for a disk-like domain, expressing it as a persistence statistic places it within the broader vineyard framework, allowing the same forecasting methodology to extend naturally to other persistence summaries. Even though we focus on the strongest feature through $M_t$, the same forecasting framework applies to any ranked vine or to aggregate summaries of several persistent vines. For domains containing essential $H_1$ classes (for example jurisdictions with enclaves or disconnected components), one would instead define $M_t$ using only finite persistence classes or use extended persistence, exactly as in previous vineyard applications such as \cite{PorterVineyards}.

We forecast our two scalar series $D_t$ and $M_t$ using standard univariate time-series methods, including seasonal autoregressive moving average (SARIMA) models \cite{hyndman2021forecasting,shumway2017time}. For each forecast, we train on completed months before a cutoff date and evaluate on the testing data (the next year). A residual bootstrap gives empirical prediction intervals by resampling training residuals and adding them to the point forecasts. We achieved a cross-validation predictive $R^2$ near 0.9 for future overdose deaths. For $M_t$, the cross-validation squared correlation with the testing data was comparable to what the $D_t$ forecast achieved.

The $M_t$ forecast can support what-if analyses that can help local health departments decide how to allocate resources. A forecast of $M_t$ estimates how severe the gap between the hardest-hit and least-hit counties is expected to become, information that can help public-health agencies anticipate periods requiring geographically targeted interventions. Public-health planners could then ask what interventions would be needed if a county-level disparity of that magnitude were to arise in different plausible locations, such as Scioto, Montgomery, Lucas, Trumbull, or another previously identified region. In each scenario, they could compare the forecasted local contrast with available treatment, harm-reduction, EMS, and outreach resources in the surrounding region, to see if the resources in the region are sufficient for various hypothetical futures. 

A forecast of $M_t$ gives information about future death spikes, but it does not identify where they will occur. We therefore also considered a simple categorical forecast for the county label associated with the strongest future $H_1$ feature. For each month, we label the most persistent finite $H_1$ class by the county containing the death simplex that kills the class, and let $Y_t$ denote this label. We then estimate a probability distribution over county labels using a categorical time-series model based on historical label frequencies weighted by recency. In our 2023 testing data, the true county label appeared in the model's top-five predicted counties in 8 of 12 months. This suggests that TDA summaries can support candidate-county forecasts for future spikes in death rates. Future researchers seeking to build on this model should include demographic, drug-market, treatment-access, reporting-delay, and policy covariates. We plan to produce improved forecasts in future work.

\section{Conclusions and Future Directions}
\label{section:Conclusion}

In this work, we applied vineyards to Ohio drug overdose data and developed a statistical framework for vineyard analysis. We introduced classical and TDA-based tests for determining when vineyards are appropriate for a spatiotemporal dataset, used vineyards to visualize the evolution of local maxima in the Ohio overdose epidemic, developed confidence bands, confidence tubes, and significance tests for vineyard features, and explored forecasting methods for both overdose deaths and topological summaries of the epidemic.

The main limitations of our work had to do with the spatial granularity of our dataset and the interpretability of the vines as they pertain to county death rates. As our data was county level, it is hard to isolate the effects of towns and cities on drug overdose deaths and limits the appropriateness of this tool for deploying resources to communities of interest or concern.  As the vines do not correspond directly to county information but rather the interplay of county death rates in the context of the whole spatial distribution of our data, our forecasting efforts cannot give us clear and specific information about when and where spikes will occur. 

Future work includes expanding on the forecasting methodology, investigating the utility of other TDA tools to analyze this dataset, determining when those methods produce statistically significant results, and adapting the statistical tests developed here to other spatiotemporal datasets.

\section*{Supplemental Materials}
\label{sec:supplemental_materials}

Additional figures and plots related to our vineyard diagrams and the spatiotemporal evolution of overdose deaths along with html files detailing all of our statistical tests can be found on our GitHub, \url{https://github.com/willeyna/OhioOverdoseVineyards}.

\acknowledgments{
The authors wish to thank Gillian Grindstaff and Mason Porter for their advice and guidance, Abigail Hickok for her code, and the anonymous reviewers for helpful comments.} 

\bibliographystyle{abbrv-doi}

\bibliography{biblio}

@article{mcmaster2025drug,
  author  = {McMaster, Ryan and Masarweh-Zawahri, Luma and Coen Flynn, Karen and Deo, Vaishali S. and Flannery, Daniel J.},
  title   = {Drug Overdose Death among Residents of Urban Census Tracts: How Granular Geographical Analyses Uncover Socioenvironmental Correlates in {C}uyahoga {C}ounty, {O}hio},
  journal = {Journal of Urban Health},
  volume  = {102},
  pages   = {445--458},
  year    = {2025},
  doi     = {10.1007/s11524-024-00939-8}
}

@book{hilbe2011negative,
  author    = {Hilbe, Joseph M.},
  title     = {Negative Binomial Regression},
  edition   = {2},
  publisher = {Cambridge University Press},
  year      = {2011},
  doi       = {10.1017/CBO9780511973420}
}

@article{yang2025topological,
  author  = {Yang, Jingjie and Fang, Heidi and Dhesi, Jagdeep and Yoon, Iris H. R. and Bull, Joshua A. and Byrne, Helen M. and Harrington, Heather A. and Grindstaff, Gillian},
  title   = {Topological Classification of Tumour-Immune Interactions and Dynamics},
  journal = {Journal of Mathematical Biology},
  volume  = {91},
  number  = {3},
  pages   = {25},
  year    = {2025},
  doi     = {10.1007/s00285-025-02253-6}
}

@article{noriega2023case,
  author  = {Noriega, Ivette and Bhullar, Manreet K. and Gilson, Thomas P. and Flannery, Daniel J. and Deo, Vaishali and Fulton, Sarah},
  title   = {A Case Study for Local Data Surveillance in Opioid Overdose Fatalities in {C}uyahoga {C}ounty, {OH} 2016--2020},
  journal = {Drug and Alcohol Dependence Reports},
  volume  = {8},
  pages   = {100187},
  year    = {2023},
  doi     = {10.1016/j.dadr.2023.100187},
  pmid    = {37711883},
  pmcid   = {PMC10498169}
}

@article{chazal_michel,
  author  = {Chazal, Fr{\'e}d{\'e}ric and Michel, Bertrand},
  title   = {An Introduction to Topological Data Analysis: Fundamental and Practical Aspects for Data Scientists},
  journal = {Frontiers in Artificial Intelligence},
  volume  = {4},
  pages   = {667963},
  year    = {2021},
  doi     = {10.3389/frai.2021.667963}
}

@book{hyndman2021forecasting,
  author = {Rob J. Hyndman and George Athanasopoulos},
  title = {Forecasting: Principles and Practice},
  edition = {3},
  year = {2021},
  note = {Online textbook}
}

@misc{choudhuri2019predicting,
  author       = {Choudhuri, Sandipan and Basu, Kaustav and Thomas, Kevin and Sen, Arunabha},
  title        = {Predicting Future Opioid Incidences Today},
  year         = {2019},
  eprint       = {1906.08891},
  archivePrefix = {arXiv},
  primaryClass = {cs.LG}
}

@article{mullen2025forecasting,
  author  = {Mullen, Aaron D. and Rock, Peter and Harris, Daniel and Slavova, Svetla and Talbert, Jeffery and Bumgardner, V. K. Cody},
  title   = {Forecasting Opioid Overdose Incidents for Rapid Actionable Data for Opioid Response in {K}entucky ({RADOR-KY})},
  journal = {Drug and Alcohol Dependence},
  volume  = {267},
  pages   = {111705},
  year    = {2025},
  doi     = {10.1016/j.drugalcdep.2024.111705}
}

@article{sumetsky2021predicting,
  author  = {Sumetsky, Natalie and Mair, Christina and Wheeler-Martin, Katherine and Cerd{\'a}, Magdalena and Waller, Lance A. and Ponicki, William R. and Gruenewald, Paul J.},
  title   = {Predicting the Future Course of Opioid Overdose Mortality: An Example From Two {US} States},
  journal = {Epidemiology},
  volume  = {32},
  number  = {1},
  pages   = {61--69},
  year    = {2021},
  doi     = {10.1097/EDE.0000000000001264},
  pmid    = {33002963},
  pmcid   = {PMC7708436}
}

@book{shumway2017time,
  author = {Robert H. Shumway and David S. Stoffer},
  title = {Time Series Analysis and Its Applications: With R Examples},
  edition = {4},
  publisher = {Springer},
  year = {2017}
}

@article{cohensteiner2007stability,
  author  = {Cohen-Steiner, David and Edelsbrunner, Herbert and Harer, John},
  title   = {Stability of Persistence Diagrams},
  journal = {Discrete \& Computational Geometry},
  volume  = {37},
  number  = {1},
  pages   = {103--120},
  year    = {2007},
  doi     = {10.1007/s00454-006-1276-5}
}

@inproceedings{ertugrul2019castnet,
  author = {Ali Mert Ertugrul and Yu-Ru Lin and Tugba Taskaya-Temizel},
  title = {{CASTNet}: Community-Attentive Spatio-Temporal Networks for Opioid Overdose Forecasting},
  booktitle = {Machine Learning and Knowledge Discovery in Databases},
  pages = {432--448},
  publisher = {Springer},
  year = {2020},
  doi = {10.1007/978-3-030-46133-1_26}
}

@article{liu2021point,
  author = {Xueying Liu and Jeremy A. Carter and Brad Ray and George Mohler},
  title = {Point Process Modeling of Drug Overdoses with Heterogeneous and Missing Data},
  journal = {The Annals of Applied Statistics},
  volume = {15},
  number = {1},
  year = {2021},
  doi = {10.1214/20-AOAS1384}
}

@article{neill2018machine,
  author  = {Daniel B. Neill and William Herlands},
  title   = {Machine Learning for Drug Overdose Surveillance},
  journal = {Journal of Technology in Human Services},
  volume  = {36},
  number  = {1},
  pages   = {8--14},
  year    = {2018},
  doi     = {10.1080/15228835.2017.1416511}
}

@article{gidea2018topological,
  author = {Marian Gidea and Yuri Katz},
  title = {Topological Data Analysis of Financial Time Series: Landscapes of Crashes},
  journal = {Physica A: Statistical Mechanics and its Applications},
  volume = {491},
  pages = {820--834},
  year = {2018}
}

@book{agresti2013,
  author = {Alan Agresti},
  title = {Categorical Data Analysis},
  edition = {3},
  publisher = {Wiley},
  year = {2013}
}

@book{anselin2023,
  author = {Luc Anselin},
  title = {An Introduction to Spatial Data Science with GeoDa},
  year = {2023}
}

@article{benjamini_hochberg_1995,
  author  = {Benjamini, Yoav and Hochberg, Yosef},
  title   = {Controlling the False Discovery Rate: A Practical and Powerful Approach to Multiple Testing},
  journal = {Journal of the Royal Statistical Society: Series B (Methodological)},
  volume  = {57},
  number  = {1},
  pages   = {289--300},
  year    = {1995},
  doi     = {10.1111/j.2517-6161.1995.tb02031.x}
}

@article{PorterVineyards,
    author = {Hickok, Abigail and Needell, Deanna and Porter, Mason A.},
    title = {Analysis of Spatial and Spatiotemporal Anomalies Using Persistent Homology: Case Studies with {COVID}-19 Data},
    journal = {SIAM Journal on Mathematics of Data Science},
    volume = {4},
    number = {3},
    pages = {1116-1144},
    year = {2022},
    doi = {10.1137/21M1435033},
    URL = {https://doi.org/10.1137/21M1435033},
    eprint = {https://doi.org/10.1137/21M1435033}
}

@Article{Cohen-Steiner2009,
author={Cohen-Steiner, David
and Edelsbrunner, Herbert
and Harer, John},
title={Extending Persistence Using {P}oincar{\'e} and {L}efschetz Duality},
journal={Foundations of Computational Mathematics},
year={2009},
month={Feb},
day={01},
volume={9},
number={1},
pages={79-103},
issn={1615-3383},
doi={10.1007/s10208-008-9027-z},
url={https://doi.org/10.1007/s10208-008-9027-z}
}

@Article{Carlsson2010,
author={Carlsson, Gunnar
and de Silva, Vin},
title={Zigzag Persistence},
journal={Foundations of Computational Mathematics},
year={2010},
month={Aug},
day={01},
volume={10},
number={4},
pages={367-405},
issn={1615-3383},
doi={10.1007/s10208-010-9066-0},
url={https://doi.org/10.1007/s10208-010-9066-0}
}

@inproceedings{VinesandVineyards,
author = {Cohen-Steiner, David and Edelsbrunner, Herbert and Morozov, Dmitriy},
title = {Vines and vineyards by updating persistence in linear time},
year = {2006},
isbn = {1595933409},
publisher = {Association for Computing Machinery},
address = {New York, NY, USA},
url = {https://doi.org/10.1145/1137856.1137877},
doi = {10.1145/1137856.1137877},
booktitle = {Proceedings of the Twenty-Second Annual Symposium on Computational Geometry},
pages = {119–126},
numpages = {8},
location = {Sedona, Arizona, USA},
series = {SCG '06}
}

@article{PHT,
    author = {Turner, Katharine and Mukherjee, Sayan and Boyer, Doug M.},
    title = {Persistent homology transform for modeling shapes and surfaces},
    journal = {Information and Inference: A Journal of the IMA},
    volume = {3},
    number = {4},
    pages = {310-344},
    year = {2014},
    month = {12},
    issn = {2049-8764},
    doi = {10.1093/imaiai/iau011},
    url = {https://doi.org/10.1093/imaiai/iau011},
    eprint = {https://academic.oup.com/imaiai/article-pdf/3/4/310/2152818/iau011.pdf},
}

@Article{Turner2024,
author={Turner, Katharine
and Robins, Vanessa
and Morgan, James},
title={The extended persistent homology transform of manifolds with boundary},
journal={Journal of Applied and Computational Topology},
year={2024},
month={Nov},
day={01},
volume={8},
number={7},
pages={2111-2154},
issn={2367-1734},
doi={10.1007/s41468-024-00175-8},
url={https://doi.org/10.1007/s41468-024-00175-8}
}

@INPROCEEDINGS{XPHTSD,
  author={Bermingham, Nicholas and Robins, Vanessa and Turner, Katharine},
  booktitle={2023 Topological Data Analysis and Visualization (TopoInVis)}, 
  title={Planar Symmetry Detection and Quantification using the Extended Persistent Homology Transform}, 
  year={2023},
  volume={},
  number={},
  pages={1-9},
  doi={10.1109/TopoInVis60193.2023.00007}}

@INPROCEEDINGS {OhioOverdoseMapper,
author = { Bermingham, Nicholas and White, David and Willey, Nathan },
booktitle = { 2025 IEEE Workshop on Topological Data Analysis and Visualization (TopoInVis) },
title = {{ Tracking the Spatiotemporal Spread of the {O}hio Overdose Epidemic with Topological Data Analysis }},
year = {2025},
volume = {},
ISSN = {},
pages = {22-31},
doi = {10.1109/TopoInVis68599.2025.00007},
url = {https://doi.ieeecomputersociety.org/10.1109/TopoInVis68599.2025.00007},
publisher = {IEEE Computer Society},
address = {Los Alamitos, CA, USA},
month =Nov}

@article{marshall2022preventing,
  title={Preventing Overdose Using Information and Data from the Environment ({PROVIDENT}): protocol for a randomized, population-based, community intervention trial},
  author={Marshall, Brandon DL and Alexander-Scott, Nicole and Yedinak, Jesse L and Hallowell, Benjamin D and Goedel, William C and Allen, Bennett and Schell, Robert C and Li, Yu and Krieger, Maxwell S and Pratty, Claire and others},
  journal={Addiction},
  volume={117},
  number={4},
  pages={1152--1162},
  year={2022},
  publisher={Wiley Online Library}
}

@misc{turner2023representingvineyardmodules,
      title={Representing Vineyard Modules}, 
      author={Katharine Turner},
      year={2023},
      eprint={2307.06020},
      archivePrefix={arXiv},
      primaryClass={math.RT},
      url={https://arxiv.org/abs/2307.06020}, 
}

@article{rosenblum2020rapidly,
  title={The rapidly changing {US} illicit drug market and the potential for an improved early warning system: evidence from {O}hio drug crime labs},
  author={Rosenblum, Daniel and Unick, Jay and Ciccarone, Daniel},
  journal={Drug and alcohol dependence},
  volume={208},
  pages={107779},
  year={2020},
  publisher={Elsevier}
}

@article{bci,
  title={A statistical analysis of drug seizures and opioid overdose deaths in {O}hio from 2014 to 2018},
  author={Ma, Lin and Tran, Lam and White, David},
  journal={Journal of Student Research},
  volume={10},
  number={1},
  year={2021}
}

@incollection{curtis2025using,
  title={Using Spatial Mixed Methods to Reveal the Geographic Nuances of Opioid Overdose Patterns in Small and Rural Towns},
  author={Curtis, Andrew and Curtis, Jacqueline and Ajayakumar, Jayakrishnan and Jefferis, Eric},
  booktitle={New Research in Crime Modeling and Mapping Using Geospatial Technologies},
  pages={211--230},
  year={2025},
  publisher={Springer}
}

@article{eck,
  title={Utilizing Supervised Machine Learning Models for Opioid Hotspot Prediction},
  author={Eck, Adam and Muradi, Aisha and Simoya, Menard},
  journal={preprint},
  year={2024},
howpublished = {\url{https://digitalcommons.oberlin.edu/researchsymp/2024/posters/5/}}
}

@article{hepler2019latent,
  title={A latent spatial factor approach for synthesizing opioid-associated deaths and treatment admissions in {O}hio counties},
  author={Hepler, Staci and McKnight, Erin and Bonny, Andrea and Kline, David},
  journal={Epidemiology},
  volume={30},
  number={3},
  pages={365--370},
  year={2019},
  publisher={LWW}
}

@book{ji2019joint,
  title={A Joint Spatio-Temporal Model of Opioid Associated Deaths and Treatment Admissions in {O}hio},
  author={Ji, Yixuan},
  year={2019},
  publisher={Wake Forest University}
}

@article{kline2021multivariate,
  title={A multivariate spatio-temporal model of the opioid epidemic in {O}hio: a factor model approach},
  author={Kline, David and Ji, Yixuan and Hepler, Staci},
  journal={Health Services and Outcomes Research Methodology},
  volume={21},
  number={1},
  pages={42--53},
  year={2021},
  publisher={Springer}
}

@article{kline2021estimating,
  title={Estimating the burden of the opioid epidemic for adults and adolescents in {O}hio counties},
  author={Kline, David and Hepler, Staci A},
  journal={Biometrics},
  volume={77},
  number={2},
  pages={765--775},
  year={2021},
  publisher={Oxford University Press}
}

@article{choi2022spatial,
  title={Spatial clustering of heroin-related overdose incidents: a case study in {C}incinnati, {O}hio},
  author={Choi, Jung Im and Lee, Jinha and Yeh, Arthur B and Lan, Qizhen and Kang, Hyojung},
  journal={BMC public health},
  volume={22},
  number={1},
  pages={1253},
  year={2022},
  publisher={Springer}
}

@article{li2019suspected,
  title={Suspected heroin-related overdoses incidents in {C}incinnati, {O}hio: A spatiotemporal analysis},
  author={Li, Zehang Richard and Xie, Evaline and Crawford, Forrest W and Warren, Joshua L and McConnell, Kathryn and Copple, J Tyler and Johnson, Tyler and Gonsalves, Gregg S},
  journal={PLoS medicine},
  volume={16},
  number={11},
  pages={e1002956},
  year={2019},
  publisher={Public Library of Science San Francisco, CA USA}
}

@article{stewart2017geospatial,
  title={Geospatial analysis of drug poisoning deaths involving heroin in the {USA}, 2000--2014},
  author={Stewart, Kathleen and Cao, Yanjia and Hsu, Margaret H and Artigiani, Eleanor and Wish, Eric},
  journal={Journal of Urban Health},
  volume={94},
  pages={572--586},
  year={2017},
  publisher={Springer}
}

@article{sudors-ohio,
  title={Unintentional Opioid Overdose Deaths in {O}hio: Insights from {S}{U}{D}{O}{R}{S} Data},
  author={White, David and Ma, Lin and Tran, Lam},
  journal={preprint},
  year={2025}
}

@article{rodriguez2023analysis,
  title={An analysis of protesting activity and trauma through mathematical and statistical models},
  author={Rodr{\'\i}guez, Nancy and White, David},
  journal={Crime Science},
  volume={12},
  number={1},
  pages={17},
  year={2023},
  publisher={Springer}
}

@article{bahid2024statistical,
  title={The statistical and dynamic modeling of the first part of the 2013-2014 {E}uromaidan protests in {U}kraine: The Revolution of Dignity and preceding times},
  author={Bahid, Yassin and Kutsenko, Olga and Rodr{\'\i}guez, Nancy and White, David},
  journal={Plos one},
  volume={19},
  number={5},
  pages={e0301639},
  year={2024},
  publisher={Public Library of Science San Francisco, CA USA}
}

@article{chen2021topological,
  title={Topological data analysis model for the spread of the coronavirus},
  author={Chen, Yiran and Voli{\'c}, Ismar},
  journal={Plos one},
  volume={16},
  number={8},
  pages={e0255584},
  year={2021},
  publisher={Public Library of Science San Francisco, CA USA}
}

@article{lo2018modeling,
  title={Modeling the spread of the {Z}ika virus using topological data analysis},
  author={Lo, Derek and Park, Briton},
  journal={PloS one},
  volume={13},
  number={2},
  pages={e0192120},
  year={2018},
  publisher={Public Library of Science San Francisco, CA USA}
}

@article{ault2022comparison,
  title={Comparison of the spread of novel coronavirus: Topological data analysis of 13 countries},
  author={Ault, Shaun V and Lu, Jia},
  journal={JIS},
  volume={6},
  number={2},
  year={2022}
}

@article{soliman2020ensemble,
  title={Ensemble forecasting of the {Z}ika space-time spread with topological data analysis},
  author={Soliman, Marwah and Lyubchich, Vyacheslav and Gel, Yulia R},
  journal={Environmetrics},
  volume={31},
  number={7},
  pages={e2629},
  year={2020},
  publisher={Wiley Online Library}
}

@article{rudkin2023spatial,
  title={Spatial Disparities in Infection Rates at the Dawn of a Pandemic: Wealthy Young Workers Mattered},
  author={Rudkin, Simon and Webber, Don J and Dlotko, Pawel},
  journal={Available at SSRN 4356837},
  year={2023}
}

@article{costatopological,
  title={A topological data analysis approach to influence-like illness},
  author={Costa, Joao Pita and {\v{S}}kraba, Primoz and Paolotti, Daniela and Mexia, Ricardo},
journal = {KDD Healthday Epidamik},
year = {2018}
}

@article{taylor2015topological,
  title={Topological data analysis of contagion maps for examining spreading processes on networks},
  author={Taylor, Dane and Klimm, Florian and Harrington, Heather A and Kram{\'a}r, Miroslav and Mischaikow, Konstantin and Porter, Mason A and Mucha, Peter J},
  journal={Nature communications},
  volume={6},
  number={1},
  pages={7723},
  year={2015},
  publisher={Nature Publishing Group UK London}
}

@Article{Ciocanel2021,
author={Ciocanel, Maria-Veronica
and Juenemann, Riley
and Dawes, Adriana T.
and McKinley, Scott A.},
title={Topological Data Analysis Approaches to Uncovering the Timing of Ring Structure Onset in Filamentous Networks},
journal={Bulletin of Mathematical Biology},
year={2021},
month={Jan},
day={16},
volume={83},
number={3},
pages={21},
issn={1522-9602},
doi={10.1007/s11538-020-00847-3},
url={https://doi.org/10.1007/s11538-020-00847-3}
}

@misc{dataohio,
  author       = {{Ohio Department of Health}},
  title        = {Mortality Data Portal},
  year         = {2024},
  howpublished = {\url{https://data.ohio.gov/wps/portal/gov/data/view/mortality}}
}

@misc{census-pop,
  author       = {{U.S. Census Bureau}},
  title        = {County Population Totals: 2010--2019. Population and Housing Unit Estimates},
  year         = {2020},
  howpublished = {\url{https://www.census.gov/data/datasets/time-series/demo/popest/2010s-counties-total.html}}
}

@article{buchanich2018effect,
  title={The effect of incomplete death certificates on estimates of unintentional opioid-related overdose deaths in the United States, 1999-2015},
  author={Buchanich, Jeanine M and Balmert, Lauren C and Williams, Karl E and Burke, Donald S},
  journal={Public Health Reports},
  volume={133},
  number={4},
  pages={423--431},
  year={2018},
  publisher={SAGE Publications Sage CA: Los Angeles, CA}
}

@article{SamuelByers2023,
  title={Detecting spatial dependence with persistent homology},
  author={Samuel Byers and Neil Pritchard and Jana Turner and Thomas Weighill},
  journal={Nonlinear Theory and Its Applications, IEICE},
  volume={14},
  number={2},
  pages={106-125},
  year={2023},
  doi={10.1587/nolta.14.106}
}

@article{pillmill,
  title = {`{T}he pill mill of {A}merica': where drugs mean there are no good choices, only less awful ones},
  author = {Arnade, Chris},
  year = {2017},
  journal = {The Guardian},
  howpublished = {\url{https://www.theguardian.com/society/2017/may/17/drugs-opiod-addiction-epidemic-portsmouth-ohio}}
}

@article{fasy,
author = {Brittany Terese Fasy and Fabrizio Lecci and Alessandro Rinaldo and Larry Wasserman and Sivaraman Balakrishnan and Aarti Singh},
title = {{Confidence sets for persistence diagrams}},
volume = {42},
journal = {The Annals of Statistics},
number = {6},
publisher = {Institute of Mathematical Statistics},
pages = {2301 -- 2339},
year = {2014},
doi = {10.1214/14-AOS1252},
URL = {https://doi.org/10.1214/14-AOS1252}
}

@article{bubenik,
author = {Bubenik, Peter},
title = {Statistical topological data analysis using persistence landscapes},
year = {2015},
issue_date = {January 2015},
publisher = {JMLR.org},
volume = {16},
number = {1},
issn = {1532-4435},
journal = {J. Mach. Learn. Res.},
month = jan,
pages = {77–102},
numpages = {26}
}

@article{chazal_bootstrap,
  author  = {Chazal, Fr{\'e}d{\'e}ric and Fasy, Brittany T. and Lecci, Fabrizio and Rinaldo, Alessandro and Singh, Aarti and Wasserman, Larry},
  title   = {On the Bootstrap for Persistence Diagrams and Landscapes},
  journal = {Modeling and Analysis of Information Systems},
  volume  = {20},
  number  = {6},
  pages   = {111--120},
  year    = {2013},
  doi     = {10.18255/1818-1015-2013-6-111-120},
  eprint  = {1311.0376},
  archivePrefix = {arXiv}
}

@article{chazal_landscapes,
  author  = {Chazal, Fr{\'e}d{\'e}ric and Fasy, Brittany T. and Lecci, Fabrizio and Rinaldo, Alessandro and Wasserman, Larry},
  title   = {Stochastic Convergence of Persistence Landscapes and Silhouettes},
  journal = {Journal of Computational Geometry},
  volume  = {6},
  number  = {2},
  pages   = {140--161},
  year    = {2015},
  doi     = {10.20382/jocg.v6i2a8}
}

@article{robinson_turner,
  title={Hypothesis testing for topological data analysis},
  author={Andrew P. Robinson and Katharine Turner},
  journal={Journal of Applied and Computational Topology},
  year={2013},
  volume={1},
  pages={241-261},
  url={https://api.semanticscholar.org/CorpusID:17167037}
}

@article{moon_lazar,
  author  = {Moon, Chul and Lazar, Nicole A.},
  title   = {Hypothesis Testing for Shapes Using Vectorized Persistence Diagrams},
  journal = {Journal of the Royal Statistical Society Series C: Applied Statistics},
  volume  = {72},
  number  = {3},
  pages   = {628--648},
  year    = {2023},
  doi     = {10.1093/jrsssc/qlad024}
}

@article{abdallah_salch,
title = {Statistical inference for persistent homology applied to simulated fMRI time series data},
author = {Hassan Abdallah and Adam Regalski and Mohammad Behzad Kang and Maria Berishaj and Nkechi Nnadi and Asadur Chowdury and Vaibhav A. Diwadkar and Andrew Salch},
journal = {Foundations of Data Science},
volume = {5},
number = {1},
pages = {1-25},
year = {2023},
issn = {},
doi = {10.3934/fods.2022014},
url = {https://www.aimsciences.org/article/id/62e247312d80b75987612297}
}

\end{document}